\documentstyle[11pt]{article}
\newcommand{\R}{{\bf R}}
 
\newcommand{\ep}{\varepsilon}
\newcommand{\rep}{r^{\varepsilon}}
\newcommand{\qep}{q^{\varepsilon}}
 
\newcommand{\Oaep}{\Omega ^{a\varepsilon}}
\newcommand{\Obep}{\Omega ^{b\varepsilon}}
\newcommand{\Oma}{\Omega^a}
\newcommand{\Omb}{\Omega ^b}
\newcommand{\Oep}{\Omega^{\varepsilon}}
\newcommand{\Jep}{J^{\varepsilon}}
\newcommand{\Taep}{T^{a\varepsilon}}
\newcommand{\Tbep}{T^{b\varepsilon}}
\newcommand{\Bbep}{B^{b\varepsilon}}
\newcommand{\Saep}{\Sigma^{a\varepsilon}}
\newcommand{\Sbep}{\Sigma^{b\varepsilon}}

\newcommand{\Hep}{H^{\varepsilon}}
\newcommand{\faep}{f^{a\varepsilon}}
\newcommand{\fbep}{f^{b\varepsilon}}

\newcommand{\haep}{h^{a\varepsilon}}
\newcommand{\hbepp}{h^{b\varepsilon}_+}
\newcommand{\hbepm}{h^{b\varepsilon}_-}
\newcommand{\oma}{\omega ^{a}}
\newcommand{\omb}{\omega ^b}
\newcommand{\repo}{r^{\varepsilon}\omega^a}
\newcommand{\kep}{k^{\varepsilon}}
\newcommand{\eaep}{e^{a\varepsilon}}
\newcommand{\ebep}{e^{b\varepsilon}}

\newcommand{\oeaep}{\overline{e}^{a\ep}}
\newcommand{\oebep}{\overline{e}^{b\ep}}
\newcommand{\oeep}{\overline{e}^{\ep}}
\newcommand{\ouep}{\overline{u}^{\varepsilon}}
\newcommand{\ouaep}{\overline{u}^{a\varepsilon}} 
\newcommand{\oubep}{\overline{u}^{b\varepsilon}}
\newcommand{\utaep}{\tilde{u}^{a\varepsilon}} 
\newcommand{\utbep}{\tilde{u}^{b\varepsilon}}
\newcommand{\etaep}{\tilde{e}^{a\varepsilon}}
\newcommand{\etbep}{\tilde{e}^{b\varepsilon}}
\newcommand{\utep}{\tilde{u}^{\varepsilon}} 
\newcommand{\etep}{\tilde{e}^{\varepsilon}}
  
\newcommand{\oza}{\overline{z}^{a}} 
\newcommand{\ozb}{\overline{z}^{b}}
\newcommand{\ou}{\overline{u}}  
\newcommand{\oua}{\overline{u}^{a}} 
\newcommand{\oub}{\overline{u}^{b}}
\newcommand{\ove}{\overline{e}}  
\newcommand{\oea}{\overline{e}^{a}} 
\newcommand{\oeb}{\overline{e}^{b}}
\newcommand{\ova}{\overline{v}^{a}} 
\newcommand{\ovb}{\overline{v}^{b}}
\newcommand{\owa}{\overline{w}^{a}} 
\newcommand{\owb}{\overline{w}^{b}}
\newcommand{\oc}{\overline{c}}
\newcommand{\ozta}{\overline{\zeta}^a}

\newcommand{\uaep}{ u^{a\varepsilon}}
\newcommand{\ubep}{ u^{b\varepsilon}}
\newcommand{\vaep}{ v^{a\varepsilon}}
\newcommand{\mepa}{ m^{\varepsilon}_\alpha}

\newcommand{\gaep}{g^{a\varepsilon}}
\newcommand{\gbep}{g^{b\varepsilon}}

\newcommand{\lep}{\lambda^\varepsilon}
\newcommand{\laep}{\lambda^{a\varepsilon}}
\newcommand{\lbep}{\lambda^{b\varepsilon}}
\newcommand{\Uep}{U^\varepsilon}
\newcommand{\Dep}{D^\varepsilon}
\newcommand{\Eep}{E^\varepsilon}
\newcommand{\eepa}{e^\varepsilon_\alpha}
\newcommand{\Eepa}{E^\varepsilon_\alpha}
\newcommand{\Depa}{D^\varepsilon_\alpha}

\newtheorem{lemma}{Lemma}
\newtheorem{theo}{Theorem}

\newtheorem{rem}{Remark}
\newtheorem{cor}{Corollary}
 
\gdef\RR{\hbox{\rm I\kern-.25em\hbox{R}}}
\gdef\NN{\hbox{\rm I\kern-.25em\hbox{N}}}
\gdef\LLI{\hbox{\rm I\kern-.25em\hbox{L}}}

\gdef\CC{\rm \hbox{C\kern-.56em\raise.4ex
         \hbox{$\scriptscriptstyle |$}\kern+0.5 em }}
\gdef\QQ{\rm \hbox{Q\kern-.35em\raise.4ex
         \hbox{$\scriptscriptstyle |$}\kern+0.5 em }}

\begin{document}
\author{Antonio GAUDIELLO, R\'egis MONNEAU, Jacqueline MOSSINO, \\ Fran\c{c}ois 
MURAT and Ali SILI}
 \date{}
\title{ {\bf Junction of elastic plates and beams \\
(Preliminary version)}} 
\maketitle
\noindent {A.G.:  
Dipartimento di Automazione, Elettromagnetismo, Ingegneria dell'Informazione e 
Matematica Industriale, Universit\`{a} di Cassino, 
Via G. Di Biasio 43, 03043 Cassino (FR), Italia; e-mail: gaudiell@unina.it} \\
\noindent {R.M: CERMICS, Ecole Nationale des Ponts et Chauss\'ees, 
6 et 8 Avenue Blaise Pascal, Cit\'e Descartes, 77455 Champs-sur-Marne Cedex 2, 
France; e-mail: monneau@cermics.enpc.fr} \\
\noindent { J.M.: C.M.L.A, Ecole Normale Sup\'erieure de Cachan, 
61 Avenue du Pr\'esident Wilson,   94235 Cachan Cedex, France; e-mail: 
mossino@cmla.ens-cachan.fr }  \\
\noindent { F.M.: Laboratoire Jacques-Louis Lions, Universit\'e Pierre et 
Marie Curie, 
Bo\^{\i}te courrier 187, 75252 Paris Cedex 05, France; e-mail: murat@ann.jussieu.fr}  \\
\noindent { A.S.: D\'epartement de Math\'ematiques, 
Universit\'e de Toulon et du Var, BP 132, 83957 La Garde Cedex, France; e-mail: 
sili@univ-tln.fr } \\ \\
{ \bf Abstract} - We consider the linearized elasticity 
system in a multidomain of $\R ^3$. This multidomain is the union of  
a horizontal plate with fixed cross section and small thickness $\varepsilon$, 
and of a vertical beam with 
fixed height and small cross section of radius $\rep$. The lateral boundary of the plate and 
the top of the beam are assumed to be clamped. When $\ep$ and $\rep$ tend to zero 
simultaneously, with $\rep \gg  \ep^2$, we identify the limit problem. 
This limit problem involves six junction conditions.\\ \\


{\bf Table of content} 

1 Introduction  \dotfill\   p. 2

2 The result \dotfill\ p. 4

\hspace{1cm} 2.1 The rescaled problem \dotfill\ p. 4

\hspace{1cm} 2.2 The setting of the limit problem \dotfill\ p. 5

\hspace{1cm} 2.3 The main result \dotfill\ p. 6

\hspace{1cm} 2.4 Back to the problem in the thin multidomain \dotfill\ p. 8

3 The derivation of the rescaled problem \dotfill\ p. 8

\hspace{1cm} 3.1 The result of the scaling \dotfill\ p. 9

\hspace{1cm} 3.2 The derivation of the scaling \dotfill\ p. 10

4 The {\em a priori} estimates and the compactness arguments \dotfill\ p. 11

\hspace{1cm} 4.1 {\em A priori} estimates \dotfill\ p. 11

\hspace{1cm} 4.2 Compactness arguments \dotfill\ p. 12

5 The limit constraints that are due to the junction \dotfill\ p. 14

\hspace{1cm} 5.1 Proof of $\oua_\alpha (0)=0$ \dotfill\ p. 14

\hspace{1cm} 5.2 Proof of $\oua_3 (x',0)\equiv\oub_3(0)$ \dotfill\ p. 14

\hspace{1cm} 5.3 Proof of $\oc(0)=0$ \dotfill\ p. 17

6 The use of convenient test functions \dotfill\ p. 20

\hspace{1cm} 6.1 The case $q=+\infty$ \dotfill\ p. 20

\hspace{1cm} 6.2 The case $q=0$ \dotfill\ p. 21

\hspace{1cm} 6.3 The case $q \in (0,+\infty)$ \dotfill\ p. 22

%7 The density arguments leading to the limit variational equation

7 Proof of stronger convergences and proof of Corollary 1 \dotfill\ p. 27

8 Appendix \dotfill\ p. 30

\hspace{1cm} 8.1 The definitions of $(v^a, w^a)$ and  
$(v^b, w^b)$ as suitable limits \dotfill\ p. 30

\hspace{1cm} 8.2 The density arguments  \dotfill\ p. 34

References \dotfill\ p. 37
 
\section{Introduction}

Let $\oma$ and $\omb$ ($a$ for "above", $b$ for "below") be two bounded 
regular domains in $\R ^{2}$. In the whole paper, the origin and axes 
are chosen so that  
\begin{equation}
\label{conddom}
\int_{\oma}x_1 \, dx_1 \, dx_2 =
\int_{\oma}x_2 \, dx_1 \, dx_2 =
\int_{\oma}x_1 x_2 \, dx_1 \, dx_2 =0 
\mbox{ and } 0 \in  \omb. 
\end{equation}
Let $\ep$ be a parameter taking values in a sequence of positive 
numbers converging to zero, and let $\rep$ be another positive parameter 
tending to zero with $\ep$. We introduce the thin multidomain $\Oep = \Oaep 
\bigcup \Jep \bigcup \Obep$, 
where $\Oaep = \repo \times (0,1)$ represents a vertical beam with 
fixed height and small cross section, $\Obep = \omb \times (-\ep, 0)$ represents a horizontal plate with small thickness and 
fixed cross section, and $\Jep=\repo \times\{0\}$ 
represents the interface at the junction between the beam and the plate.

In this thin multidomain, we consider the displacement 
$\overline{U}^\ep$, solution of the 
three-dimensional linearized elasticity system:
\begin{eqnarray}  \label{varpetit}
\overline{U}^\ep \in  Y^\ep\mbox{ and } 
\forall U \in  Y^\ep, 
\int_{\Oep} [A^\ep e(\overline{U}^\ep), e(U)]\, dx 
&=&
\int_{\Oep}F^\ep. U\, dx +
\int_{\Oep} [G^\ep, e(U)]\, dx + 
\nonumber \\ &+&
\int_{\Saep \bigcup \Tbep \bigcup \Bbep}H^\ep. U\, d\sigma,
\end{eqnarray}
where 

$\bullet$ $ Y^\ep=\{U \in (H^1(\Oep))^3,\, U=0  \mbox{ on } \Taep=\repo \times \{ 1\} 
 \mbox{ and on }\Sbep = \partial \omb \times (-\ep,0)\}$;

$\bullet$ $A^\ep =A^\ep(x)=\left \{ \begin{array}{c}
  A^a, \mbox{ if } x \in \Oaep, \\
\kep  A^b, \mbox{ if } x \in \Obep, \end{array}\right. \\$
with $\kep $ a positive parameter depending on $\ep$ and $A^a$, $A^b$ 
tensors with constant coefficients $A^a_{ijkl}$ and $A^b_{ijkl}$, 
$i, j, k, l \in \{1,2,3\}$, satisfying the usual symmetry and coercivity 
conditions: 
\[
A^a_{ijkl}=A^a_{jikl}=A^a_{ijlk}, \quad 
A^b_{ijkl}=A^b_{jikl}=A^b_{ijlk},
\]
\[%\begin{equation}\label{coer} 
\exists C>0, \; \forall \xi \in \R^{3 \times 3}_s, 
[A^a \xi, \xi] \geq  C |\xi|^2, \; [A^b \xi, \xi] \geq  C |\xi|^2,
\]%\end{equation}
where $\R^{3 \times 3}_s$ denotes the set of symmetric $3\times 3$-matrices, 
$(A^a \xi)_{ij}=  \sum_{kl}A^a_{ijkl} \xi_{kl}$ (e.g.), 
the scalar product $[.,.]$ in $\R^{3 \times 3}_s$ is defined by 
$[\eta,\xi]= \sum_{ij}\eta_{ij}\xi_{ij}$ 
and $|.|$ is the associated norm;\\

$\bullet$ $e_{ij}(U) = \frac{1}{2} \left( 
\frac{\partial U_i}{\partial x_j} +\frac{\partial U_j}{\partial x_i} \right)$;\\

$\bullet$ $F^\ep \in (L^2 (\Oep))^{3}$; the euclidian scalar product in 
$\R^{3}$ is denoted by a dot;\\

$\bullet$ $G^\ep \in (L^2 (\Oep))^{3 \times 3}$;\\ 
 
$\bullet$ $H^\ep \in (L^2(\Saep \bigcup \Tbep \bigcup \Bbep ))^{3}$, 
where $\Saep$ denotes the lateral boundary of the beam, $\Tbep$ and 
$\Bbep$ are respectively  the top and the bottom of the plate: 
$$\Saep  =  \rep \partial \oma \times (0,1), \quad
\Tbep = (\omb \setminus \repo ) \times \{0\}, \quad
\Bbep = \omb \times \{-\ep\}.$$

The constraint "$U=0$" in the definition of $ Y^\ep$ means that the 
multistructure is clamped on the top $\Taep$ of the beam and on the 
lateral boundary $\Sbep$ of the plate. The case $\kep$ tending 
to zero or infinity 
corresponds to very different materials in $\Oaep$ and $\Obep$. (Note that 
breaking the symmetry between $\Oaep$ and $\Obep$ is not restrictive.) In the right 
hand side of (\ref{varpetit}), the second term is written in divergence form 
like in \cite{MuSi1}, \cite{MuSi2} and \cite{GaMoMoMuSi}. It is well known 
that, by means of the Green formula, this second term can contribute to 
the other ones, giving possibly less regular (not necessarily $L^2$) volume and surface 
source terms. For convenience of the reader, we have chosen to write the three 
integrals: one recovers the classical formulation by setting $G^\ep=0$, but 
the simplest case corresponds to $F^\ep=0, H^\ep=0$ and $G^\ep \neq 0$. This case was considered in 
the short preliminary version \cite{GaMoMoMuSi}.

Problem (\ref{varpetit}) admits a unique solution 
$\overline{U}^\ep$ (see  e.g. \cite{OlShYo}). The aim of this paper is to 
describe the limit behaviour of the displacement $\overline{U}^\ep$, 
as $\ep$ tends to zero. We prove that this behaviour depends on the 
limit of the sequence 
 \[
\qep=  \kep \frac{\ep^3}{(\rep)^2}.
 \]
When $\kep  \ep^3$ and $(\rep)^2$ have same order (i.e. $\qep \rightarrow q \in 
(0,\infty)$), the limit problem (after suitable rescaling) is coupled 
between a two-dimensional plate and a one-dimensional beam, with six 
junction conditions. If $\kep  \ep^3 \gg (\rep)^2$, the multistructure 
has the limit behaviour of a thin rigid plate and a thin elastic beam which 
are independent of each other, the beam being clamped at both ends; 
on the contrary, if $\kep  \ep^3 \ll   (\rep)^2$, the structure behaves as a 
thin rigid beam and a thin elastic plate which are independent of 
each other, the plate being clamped on its contour and fixed vertically 
at the junction.
 
The reader is referred to \cite{AcBuPe}, \cite{AnBaPe}, \cite{Ca}, \cite{Ci2}, \cite{CiDe}, 
\cite{DaGr}, \cite{FrJaMu}, \cite{Le2}, \cite{LeRa1}, \cite{LeRa2}, \cite{MuSi1}, 
\cite{MuSi2}, \cite{Pe}, \cite{TrVi}, for the asymptotic behaviour of plates and 
beams. 
Junction problems are considered in \cite{Ci1}, \cite{CiSa}, \cite{GaGuLeMo2}, 
\cite{GaGuLeMo3}, \cite{Gr}, \cite{KoMaMo}, \cite{Le1}. The present work is a natural follow 
up of \cite{MuSi1}, \cite{MuSi2}, which deal with reduction of dimension for 
elastic thin cylinders, and \cite{GaGuLeMo2}, \cite{GaGuLeMo3}, which deal 
with the diffusion 
equation in the thin multistructure considered in this paper. Our results were 
announced in a short note \cite{GaMoMoMuSi}.
 
\section{The result}

\subsection{The rescaled problem}

In the sequel, the indexes $\alpha$ and $\beta$ take values in the set 
$\{1,2\}$. Moreover, $x=(x',x_3)$ denotes the generic point in $\R^3$.

Let $\Oma = \oma \times (0,1)$, $\Omb = \omb \times (-1,0)$, 
$T^a=\oma \times \{1\}$, 
%$T^b=\omb \times \{0\}$, $B^b=\omb \times \{-1\}$, 
$\Sigma^a= \partial \oma \times (0,1)$ and $\Sigma^b = \partial \omb \times 
(-1,0)$. The asymptotic behaviour of $\overline{U}^\ep$ can be described by 
using a convenient rescaling. (The reader is referred to Section 3.1 for details.)
This rescaling maps the space $ Y^\ep$ onto the space ${\cal Y} ^\ep $ 
defined by:
\begin{equation}\label{defy}\begin{array}{l}
{\cal Y} ^\ep=\{ u= (u^a,u^b) \in (H^1 (\Oma))^3 \times (H^1 (\Omb))^3, \;
u^a =0 \mbox{ on } T^a, u^b =0 \mbox{ on }\Sigma^b,\\
\quad \quad \quad\mbox{ for a.e. }x' \in \oma, \;
u^a_\alpha (x',0) = \ep \rep u^b_\alpha (\rep x',0) \mbox{ and }\;
u^a_3(x',0) =  u^b_3(\rep x',0) \}.\end{array} 
\end{equation}
In particular, we denote by $\ouep=(\ouaep,\oubep)$ the rescaling of 
the solution $\overline{U}^\ep$ of Problem (\ref{varpetit}). We set
\begin{equation}\label{oeep}
 e^{a\ep}(u^a) = \left(
\begin{array}{cc} \frac{1}{(\rep)^2} e_{\alpha \beta}(u^a) 
& \frac{1}{\rep}e_{\alpha  3}(u^a ) \\ \\
\frac{1}{\rep}e_{3 \alpha}( u^a  )
& e_{33}(u^a )\end{array}\right), \; \;
 e^{b\ep}(u^b) = \left(
\begin{array}{cc}  e_{\alpha \beta}(u^b) 
& \frac{1}{\ep} e_{\alpha  3}(u^b ) \\ \\
\frac{1}{\ep} e_{3 \alpha}(u^b )
& \frac{1}{\ep^2}e_{33}( u^b )\end{array}\right).
\end{equation}
Then $\ouep$ 
is the unique solution of the following problem:

\begin{equation}\label{vargrand} \begin{array}{l}
\ouep \in {\cal Y} ^\ep\mbox{ and } \forall u \in  {\cal Y} ^\ep, \\ \\
\quad \int_{\Oma} [A^a e^{a \ep}(\ouaep) , e^{a\ep}(u^a) ]\, dx + \qep 
\int_{\Omb} [A^b e^{b \ep}(\oubep) , e^{b\ep}(u^b) ]\, dx= \\ \\
\quad \quad =\int_{\Oma}\faep . u^a\, dx +\int_{\Omb}\fbep . u^b\, dx +
 \int_{\Oma}[\gaep , e^{a\ep}(u^a)]\, dx +
\int_{\Omb}[\gbep , e^{b\ep}(u^b)]\, dx +\\ \\
\quad  \quad \quad\int_{\Sigma^a}\haep . u^a\, d\sigma  + 
\int_{\omb} \left( \hbepp . u^b_{ \;| x_3 =0}  
+ \hbepm . u^b_{\;  | x_3 =-1} \right)\,dx',
\end{array}
\end{equation}
where $\qep$ is defined by 
\begin{equation}\label{qep}
\qep=  \kep \frac{\ep^3}{(\rep)^2}
\end{equation}
and where the source terms are suitable transforms of 
$(F^\ep, G^\ep, H^\ep)$ (see Section 3.1).
 
\subsection{The setting of the limit problem}

For the definition of the limit problem, 
%in a way similar to \cite{MuSi1} and \cite{MuSi2}, 
we introduce the following functional spaces: 
\[
%\begin{equation}\label{calUa}
\begin{array}{ll}

{\cal U}^a= 
&\{u^a \in (H^2_0(0,1))^2 \times H^1(\Oma),   
%\;
\exists \zeta^a \in H^1(0,1), 
%\;
\zeta^a (1)=0, 
%\quad 
u^a_3 = \zeta^a -
x_{1}\frac{d\,u^a_1}{d\,x_3} -
x_{2}\frac{d\,u^a_2}{d\,x_3}
\},\\
\\
 
%\end{equation}
%\begin{equation}\label{calVa}

{\cal V}^a= 
&\{
 v^a \in (H^1(\Oma))^2 \times L^2(0,1;H^1 (\oma)), \; 
\exists c \in H^1_0(0,1),
%\\  
%&\ \   
v^a_1=-c\,x_2, v^a_2=c\,x_1, 
%\quad 
\\&  \ \
\mbox{ for a.e. } 
x_3 \in (0,1), \; 
 \int_{\oma}v^a_3 (x',x_3)\, dx'=0 
\},\\
\\

%\end{equation}

%\begin{equation}\label{calWa}

{\cal W}^a= 
&\{
w^a \in (L^2(0,1;H^1 (\oma)))^2 \times \{0\}, \\
&  \ \  \mbox{for a.e. } x_3 \in (0,1),  \;
\int_{\oma} w^a_\alpha \, dx'  = 
\int_{\oma} (x_1 w_2^a -x_2 w_1^a)\, dx' =0
\},\\
\\

%\end{equation}

%\begin{equation}\label{calUb}
{\cal U}^b= 
&\{u^b \in (H^1(\Omb))^2 \times H^2_0(\omb), \;
\exists \zeta^b_\alpha \in H^1_0(\omb ), \quad 
u^b_\alpha = \zeta^b_\alpha -
x_3\frac{\partial u^b_3}{\partial x_\alpha}
\},\\
\\
 
%\end{equation}

%\begin{equation}\label{calVb}
{\cal V}^b= 
&\{v^b \in (L^2(\omb ;H^1 (-1,0)))^2  \times  \{0\}, \quad 
 \mbox{ for a.e. } x' \in \omb,  \; 
\int_{-1}^{0}v^b_\alpha (x',x_3)\, dx_3=0  
\},\\
\\

%\end{equation}

%\begin{equation}\label{calWb}
{\cal W}^b= 
&\{w^b \in (\{0\})^2  \times L^2(\omb ;H^1 (-1,0)) , \quad 
 \mbox{ for a.e. } x' \in  \omb,  \;
\int_{-1}^{0}  w^b_3 (x',x_3)\, dx_3  =0
\},
\end{array}
\]
%\end{equation}

\[
%\begin{equation}\label{calZab}
{\cal Z}^a = {\cal U}^a \times {\cal V}^a \times {\cal W}^a, \;  \;
{\cal Z}^b = {\cal U}^b \times {\cal V}^b \times {\cal W}^b.
\]

Without loss of generality, we assume that $\qep$ defined by (\ref{qep}) 
satisfies
\begin{equation}\label{convq}
\qep \rightarrow q \in [0, \infty].
\end{equation} 
According to the value of $q$, the functional space for the limit 
problem is the following one:
\[%\begin{equation}\label{calZ}
\begin{array}{llll}
{\cal Z}   &  =  \{ z= (z^a,z^b) = ((u^a,v^a,w^a),(u^b,v^b,w^b)) \in
                           {\cal Z}^a \times {\cal Z}^b, &\\
           &  \ \ \ \ \ \mbox{for a.e. } x' \in \oma,  \; u^a_3 (x',0) =
                           u^b_3 (0) \},  &\; \mbox{ if } q \in (0, +\infty),\\
&&\\
{\cal Z}_{\infty} & = \{  z^a =  (u^a,v^a,w^a) \in 
{\cal Z}^a, \quad  
\mbox{ for a.e. } x' \in \oma,  \;
u^a_3 (x',0) = 0 \}, &\; \mbox{ if } q = +\infty,\\
&&\\
{\cal Z}_0 &= \{  z^b =  (u^b,v^b,w^b) \in 
{\cal Z}^b,  
   \;
 u^b_3 (0)=0\},&\; \mbox{ if } q = 0.
\end{array} 
\]%\end{equation}

Let us remark that ${\cal U}^a$ (resp. ${\cal U}^b$) 
is a Bernouilli-Navier (resp. Kirchhoff-Love) space of displacements. 
Less classical spaces are ${\cal V}^a$, ${\cal W}^a$, ${\cal V}^b$, 
${\cal W}^b$, which are 
introduced in a way similar to \cite{MuSi1} and \cite{MuSi2} 
(see also Appendix, Section 8.1). As for the 
boundary conditions, some of them are due to the clamping. 
These are more or less standard ones:
$$
u^a_\alpha(1)=\frac{du^a_\alpha}{dx_3}(1)=c(1)=0,\quad u^b_3=0 \mbox{ and }
\frac{\partial u^b_3}{d\nu} =0 \mbox{ on }\partial \omb.
$$
In contrast with the other requirements, the six conditions 
$$ \displaystyle{
u^a_{\alpha}(0)=\frac{d\,u^a_{\alpha}}{d\,x_3}(0) = c(0) =0}, \quad
u^a_3 (x',0) = u^b_3 (0) \; \mbox{(respectively } u^a_3 (x',0) = 0 \mbox{ or }
u^b_3 (0) =0),$$ 
which appear in the definition of the above spaces 
${\cal U}^a$, ${\cal V}^a$ 
and ${\cal Z}$ (respectively ${\cal Z}_\infty$ or ${\cal Z}_0$), are specific 
to the junction between the beam and the plate. Note also that, 
in view of the definition of ${\cal U}^a$, the condition $u^a_3 (x',0) = u^b_3 
(0)$ (respectively $u^a_3 (x',0) = 0$) reduces to $\zeta^a(0)=
u^b_3(0)$ (respectively $\zeta^a(0)=0$).\\
 
We finally introduce, for $z^a=(u^a,v^a,w^a)$ in ${\cal Z}^a$ and 
$z^b=(u^b,v^b,w^b)$ in ${\cal Z}^b$:
\begin{equation}\label{eab}
 e^{a}(z^a) = \left(
\begin{array}{cc}  e_{\alpha \beta}(w^a) 
&  e_{\alpha  3}(v^a ) \\ \\
 e_{3 \alpha}(v^a  )
& e_{33}(u^a)\end{array}\right), \; \;
 e^{b}(z^b) = \left(
\begin{array}{cc}  e_{\alpha \beta}(u^b) 
&  e_{\alpha  3}(v^b ) \\ \\
 e_{3 \alpha}(v^b  )
& e_{33}(w^b)\end{array}\right).
\end{equation}

\subsection{The main result} 

%Let $f^{a\ep}$, $f^{b\ep}$, $g^{a\ep}$, $g^{b\ep}$, $h^{a\ep}$, 
%$\hbepp$, $\hbepm$ be defined by (\ref{fep}), (\ref{gep}), (\ref{hep}). 
We describe the limit behaviour of Problem (\ref{vargrand}), 
as $\ep$ tends to zero. In the sequel, we assume that
\begin{equation}\label{convfa}
\faep \rightharpoonup f^a \mbox{ weakly in } (L^2(\Oma ))^{3}, 
\end{equation}
\begin{equation}\label{convfb}
\fbep \rightharpoonup f^b \mbox{ weakly in } (L^2(\Omb ))^{3}, 
\end{equation}
\begin{equation}\label{convga}
\gaep \rightharpoonup g^a \mbox{ weakly in } (L^2(\Oma ))^{3 \times 3}, 
\end{equation}
\begin{equation}\label{convgb}
\gbep \rightharpoonup g^b \mbox{ weakly in } (L^2(\Omb ))^{3\times 3},
\end{equation}
\begin{equation}\label{convha}
\haep \rightharpoonup h^a \mbox{ weakly in } (L^2(\Sigma^a))^{3}, 
\end{equation}
\begin{equation}\label{convhb} 
\hbepp\rightharpoonup h^b_+ \mbox{ and }\hbepm\rightharpoonup h^b_-
\mbox{ weakly in } (L^2(\omb))^{3},
\end{equation}
which is not restrictive, as proved in Remark 2 hereafter.
 
Our main result is the following one: 
 
\begin{theo} Assume that $\frac{\rep}{\ep^2} \rightarrow + \infty$ and that 
(\ref{convq}), (\ref{convfa}) to (\ref{convhb}) hold true. 
Then, with  $\eaep$, $\ebep$ defined in (\ref{oeep}) and $e^a$, $e^b$  
defined in (\ref{eab}),\\

1) If $q \in (0, +\infty)$, there exists $\overline{z}= 
(\oza,\ozb)=
((\overline{u}^a,\overline{v}^a,\overline{w}^a),
(\overline{u}^b,\overline{v}^b,\overline{w}^b)) \in {\cal Z}$, such that
\begin{equation}
\label{convouep1} 
(\ouaep,\oubep) \rightharpoonup (\oua, \oub) \mbox{ weakly in } (H^1(\Oma))^3 
\times (H^1(\Omb))^3, 
\end{equation}
\begin{equation}\label{convoeep1}  
(\eaep (\ouaep), \ebep (\oubep)) \rightharpoonup ( e^a (\oza),e^b (\ozb)) \mbox{ weakly in }
(L^2(\Oma))^{3 \times 3}\times (L^2(\Omb))^{3 \times 3},
%\end{array}
\end{equation}
and $\overline{z}$ is the unique solution of the following problem:
\begin{equation}\label{vargrand1} \begin{array}{l}
\overline{z}  \in  {\cal Z} \mbox{ and } \forall z \in {\cal Z}, \\ \\
\quad \int_{\Oma} [A^a e^a(\oza), e^a(z^a)]\, dx+ 
 q \int_{\Omb} [A^b e^b(\ozb), e^b(z^b)]\, dx= \\ \\
\quad \quad 
=\int_{\Oma}f^a . u^a\, dx +\int_{\Omb}f^b . u^b\, dx +
\int_{\Oma} [ g^a ,e^a(z^a) ]\, dx  +
\int_{\Omb} [ g^b, e^b(z^b)]\, dx  +\\ \\
\quad  \quad \quad +\int_{\Sigma^a}h^a . u^a\, d\sigma + 
\int_{\omb} \left( h^b_+ . u^b_{ \;| x_3 =0}  
+ h ^b_- . u^b_{\;  | x_3 =-1} \right)\,dx'.
\end{array}
\end{equation}
 
Moreover, if the convergences in (\ref{convga}), (\ref{convgb}) are strong, 
then 
%$(g^a,g^b) \neq (0,0)$ 
the convergences in (\ref{convouep1}) and (\ref{convoeep1}) are strong. 
\\
 
2) If $ q=+\infty$, there exists $\oza= 
(\overline{u}^a,\overline{v}^a,\overline{w}^a)
  \in {\cal Z}_\infty$, such that
\begin{equation} \label{convouep2}
\ouaep \rightharpoonup \oua \mbox{ weakly in } (H^1(\Oma))^3,\; \;
\oubep \rightarrow 0 \mbox{ strongly in } (H^1(\Omb))^3,
\end{equation}
\begin{equation} \label{convoeep2} 
\eaep (\ouaep) \rightharpoonup e^a (\oza) \mbox{ weakly in }(L^2(\Oma))^{3 \times 3}, \;\;
 \ebep (\oubep)\rightarrow 0 \mbox{ strongly in }(L^2(\Omb))^{3 \times 3},
\end{equation}
and $\overline{z}^a$ is the unique solution of the following problem:
\begin{equation}\label{vargrand2} \begin{array}{l}
\oza \in {\cal Z}_{\infty} \mbox{ and } \forall z^a \in {\cal Z}_{\infty}, 
\\ \\  
\quad \int_{\Oma} [A^a e^a(\oza), e^a(z^a)]\, dx 
=\int_{\Oma}f^a . u^a\, dx  + 
 \int_{\Oma} [ g^a ,e^a(z^a) ]\, dx  +
 \int_{\Sigma^a}h^a . u^a\, d\sigma.
\end{array}
\end{equation}
Moreover, if the convergence in (\ref{convga}) is strong,
then
\begin{equation}\label{convu2}
\ouaep \rightarrow \oua \mbox{ strongly in } (H^1(\Oma))^3,
\end{equation}
\begin{equation}\label{conve2} 
\eaep (\ouaep) \rightarrow e^a (\oza) \mbox{ strongly in }(L^2(\Oma))^{3 \times 3}, \; \;
\sqrt{\qep} \; \ebep (\oubep) \rightarrow 0 \mbox{ strongly in }(L^2(\Omb))^{3 \times 3}.
\end{equation}

3) If $q=0$, there exists $\ozb= 
(\overline{u}^b,\overline{v}^b,\overline{w}^b)
  \in {\cal Z}_0$, such that
\begin{equation}\label{convouep3} 
\qep \, \ouaep \rightarrow 0 \mbox{ strongly in } (H^1(\Oma))^3,\;
\qep \, \oubep \rightharpoonup \oub \mbox{ weakly in } (H^1(\Omb))^3,  
\end{equation}
\begin{equation}\label{convoeep3}  
\qep \eaep (\ouaep)  \rightarrow 0 \mbox{ strongly in }(L^2(\Oma))^{3 \times 3}, \;
\qep  \ebep (\oubep)\rightharpoonup e^b (\ozb) \mbox{ weakly in }(L^2(\Omb))^{3 
\times 3}, 
\end{equation}
and $\overline{z}^b$ is the unique solution of the following problem:
\begin{equation}\label{vargrand3} \begin{array}{l}
\ozb \in {\cal Z}_0\mbox{ and } \forall z^b \in {\cal Z}_0,  \\ \\
\quad  \int_{\Omb} [A^b e^b(\ozb), e^b(z^b)]\, dx    
 =\int_{\Omb}f^b . u^b\, dx   + 
 \int_{\Omb} [ g^b ,e^b(z^b) ]\, dx  + \\ \\
\qquad \qquad+\int_{\omb} \left( h^b_+ . u^b_{ \;| x_3 =0}  
+ h ^b_- . u^b_{\;  | x_3 =-1} \right)\,dx'. 
\end{array}
\end{equation}
Moreover, if the convergence in (\ref{convgb}) is strong,  
then
\begin{equation}\label{convu3}
\qep \oubep \rightarrow \oub \mbox{ strongly in } (H^1(\Omb))^3,
\end{equation}
\begin{equation}\label{conve3}
\sqrt{\qep} \; \eaep (\ouaep) \rightarrow 0 \mbox{ strongly in }
(L^2(\Omb))^{3 \times 3}, \; \;
\qep \ebep (\oubep) \rightarrow e^b (\ozb) \mbox{ strongly in }
(L^2(\Omb))^{3 \times 3}. 
\end{equation} 
\end{theo}

\begin{rem} The condition $\frac{\rep}{\ep^2} \rightarrow + \infty$ is only 
used to prove that $\oua_3(x',0)\equiv \oub_3(0)$ and $\oc (0)=0$, via a 
convenient Sobolev imbedding theorem, as regards the second equality. At this stage of 
our understanding, we do not know if it is just a technical condition or not. 
\end{rem} 
 
\begin{rem} In the Appendix, Section 8.1, we prove that the functions 
$\ova$ and $\owa$ (resp. $\ovb$ and $\owb$) which appear in the limit problem 
are the limits of suitable expressions of $\ouaep$ (resp. 
$\oubep$).
\end{rem} 
 
\subsection{Back to the problem in the thin multidomain} 
 
As far as the asymptotic behaviour of the "energy" of the solution of 
Problem (\ref{varpetit}) in the thin 
multidomain  is concerned, we define the following renormalized energy by:
\begin{equation}\label{ren.en}
{\cal E}^\ep : = \left(\frac{\lep}{\rep}\right)^2
\int_{\Oep} [A^\ep e(\overline{U}^\ep), e(\overline{U}^\ep)]\, dx,
\end{equation} 
where $\lep$ can be made explicit in terms of $\ep, \rep, F^\ep, G^\ep, H^\ep$ 
(see (\ref{lep}) in Section 3.1); we also have 
\[{\cal E}^\ep =\int_{\Oma} [A^a \eaep(\ouaep), \eaep(\ouaep)]\, dx + 
\qep \int_{\Omb} [A^b \ebep(\oubep) , \ebep(\oubep)]\, dx 
\] 
and from Theorem 1 we deduce the following Corollary:

\begin{cor} Assume that $\frac{\rep}{\ep^2} \rightarrow + \infty$ and that 
(\ref{convq}), (\ref{convfa}) to (\ref{convhb}) hold true. \\ 

1) If $q \in (0, +\infty)$ and the convergences 
in (\ref{convga}), (\ref{convgb}) are strong, then 
\[  
{\cal E}^\ep \rightarrow {\cal E}=\int_{\Oma} [A^a e^a(\oza),e^a(\oza) ]\, dx + 
 q \int_{\Omb} [A^b e^b(\ozb), e^b(\ozb) ]\, dx.
\] 

2) If $ q=+\infty$ and the convergence in (\ref{convga}) 
is strong, then
\[ 
{\cal E}^\ep \rightarrow  {\cal E}_{\infty}=
\int_{\Oma} [A^a e^a(\oza),e^a(\oza) ]\, dx. 
\] 

3) If $q=0$ and the convergence in (\ref{convgb}) is strong, 
then
\[
\qep {\cal E}^\ep  \rightarrow {\cal E}_0 =
\int_{\Omb} [A^b e^b(\ozb),e^b(\ozb) ]\, dx.  
\]
\end{cor}
\vspace{0.5cm}
 
{\bf The remaining part of the paper is devoted to the proofs of Theorem 1 and 
Corollary 1.}
 
%\section{The proof}
%\sub
\section{The derivation of the rescaled problem}
%The computations leading to the rescaled problem are easy. For briefness, we 
%leave them to the reader. 
Let us emphasize that we perform different scalings 
for the respective restrictions of $U \in Y^\ep$ to the respective subdomains 
$\Oaep$ and $\Obep$, in order to get convenient transmission conditions 
for their transforms $u^a$ and $u^b$. We mean that, with the transmission conditions 
appearing in the definition (\ref{defy}) of ${\cal Y}^\ep$, namely 
\begin{equation}\label{tc}
\mbox{ for a.e. }x' \in \oma, \;
u^a_\alpha (x',0) = \ep \rep u^b_\alpha (\rep x',0) \mbox{ and }
u^a_3(x',0) =  u^b_3(\rep x',0),
\end{equation} 
we are abble to derive the junction conditions, for the limit problem. (The derivation of the limit junction conditions seems 
to be delicate otherwise.) Moreover this is the scaling for which the coupling is 
maximum at the limit, at least for the third component of the 
displacement.

%\sub
\subsection{The result of the scaling}
In this subsection, we give the explicit expressions of the source terms and 
the solution of the rescaled problem (\ref{vargrand}), as functions of the 
corresponding terms of the initial problem (\ref{varpetit}). An explanation 
is given in Subsection 3.2.

On the first hand, assuming that 
$(F^\ep, G^\ep, H^\ep) \neq (0,0,0)$ (otherwise the problem is trivial), 
we define $\lep$ by: 
\begin{equation}\label{lep}
\begin{array}{l}
\frac{1}{(\rep)^2} \sum_{\alpha=1}^{2}\|F^\ep_{\alpha}\| ^2_{L^2(\Oaep)} +
\|F^\ep_{3}\| ^2_{L^2(\Oaep)} +
\frac{\ep^3}{(\rep)^2} \sum_{\alpha=1}^{2}\|F^\ep_{\alpha}\| ^2_{L^2(\Obep)} +
\frac{\ep}{(\rep)^2}\|F^\ep_{3}\| ^2_{L^2(\Obep)} + \\
+
\|G^\ep\| ^2_{(L^2(\Oaep))^{3\times 3}} + \frac{\ep^3}{(\rep)^2}
\|G^\ep\| ^2_{(L^2(\Obep))^{3\times 3}}  
+\frac{1}{(\rep)^3}\sum_{\alpha=1}^{2}\|H^\ep_{\alpha}\| ^2_{L^2(\Saep)}
+\frac{1}{\rep}\|H^\ep_{3}\| ^2_{L^2(\Saep)} +\\
+\frac{\ep^2}{(\rep)^2}\sum_{\alpha=1}^{2}
\|H^\ep_{\alpha}\| ^2_{L^2(\Tbep \bigcup \Bbep)}
+\frac{1}{(\rep)^2}\|H^\ep_{3}\| ^2_{L^2(\Tbep  \bigcup \Bbep )}=
\left(\frac{\rep}{\lep}\right)^2;
\end{array}
\end{equation}
then we set
\begin{equation}\label{fep}
\begin{array}{lll}
%\begin{eqnarray} \label{fep} 
\faep_\alpha (x) = \frac{\lep}{\rep}F^\ep_\alpha ( \rep x',x_3) ,  \; \; & &
\faep _3 (x) =     \lep\, F^\ep_3
(\rep x',x_3), \quad \mbox{ for }x\in \Oma;    \\ \\
\fbep_\alpha (x)  = \lep   \frac{\ep^2}{(\rep)^2}\, F^\ep_\alpha 
(  x', \ep x_3), \;  \;  & &
\fbep _3 (x) = \lep\frac{\ep}{(\rep)^2}  \,  F^\ep_ 3
(  x', \ep x_3),\quad  \mbox{ for }  x\in \Omb; 
\end{array}
\end{equation}

\begin{equation} \label{gep}
\gaep(x)= \lep G^\ep 
(\rep x',x_3),\;  \mbox{ for }  x\in \Oma; \; \; \; \;
\gbep(x)=\lep\frac{\ep^2}{(\rep)^2} G^\ep  (x', \ep x_3), \;\mbox{ for }  x\in \Omb;
\end{equation}

\begin{equation}\label{hep} 
\begin{array}{lll}
\haep_\alpha (x) = \lep\frac{1}{(\rep)^2}H^\ep_\alpha ( \rep x',x_3), 
& \haep _3 (x) = \lep \frac{1}{\rep} H^\ep_3(\rep x' ,x_3), 
& \mbox{for } x \in  \Sigma^a;
 \\ \\
%\hbepp \mbox{ and }\hbepm : \omb \rightarrow \R \mbox{ are given by: }\\ \\
{\hbepp}_\alpha(x') = {\hbepp}_3(x')=0, \; & &\mbox{for }x' \in \repo, \\ \\
{\hbepp} _\alpha (x')  = \lep   \frac{\ep}{(\rep)^2}\, H^\ep_\alpha 
(  x', 0), 
& {\hbepp} _3 (x') = \lep\frac{1}{(\rep)^2}  \,  H^\ep_ 3(  x', 0), 
& \mbox{for }  x'\in \omb \setminus \repo, \\ \\
{\hbepm} _\alpha (x')  = \lep   \frac{\ep}{(\rep)^2}\, H^\ep_\alpha 
(  x', -\ep), 
& {\hbepm} _3 (x') =\lep \frac{1}{(\rep)^2}  \,  H^\ep_ 3
(  x',-\ep ), 
& \mbox{for }  x'\in \omb.
\end{array}
\end{equation}
(Note that $\hbepp=0$ on $\rep \oma$, since there is no contribution of $\Hep$ 
on $\Jep$.)

On the other hand, for any $U \in  Y^\ep$, we define the rescaled function 
$u=(u^a,u^b)$ by:
\begin{equation} \label{ouaep}
u^a _\alpha (x) = \lep   \rep \, U_\alpha 
(\rep x',x_3),  \; \; 
  u^a_3 (x) =     \lep\, U_3
(\rep x',x_3),  \mbox{ for }x\in \Oma,    
\end{equation} 
\begin{equation} \label{oubep}
 u^b_\alpha (x) = \lep   \frac{1}{\ep}\, U_\alpha 
(  x', \ep x_3), \;  \;  
u^b _3 (x) =   \lep  \, U_3
(  x', \ep x_3), \mbox{ for }  x\in \Omb. 
\end{equation}

\begin{rem}Let us remark that the rescaled souce terms are bounded, but not strongly converging to zero, since, by definition of $\lep$ 
(see (\ref{lep})) and by (\ref{fep}) to (\ref{hep}),
\[
\begin{array}{l}
\|\faep\|^2_{(L^2(\Oma))^{3}} + \|\fbep\|^2_{(L^2(\Omb))^{3}} + 
\|\gaep\|^2_{(L^2(\Oma))^{3 \times 3}} 
+\|\gbep\|^2_{(L^2(\Omb))^{3 \times 3}} +\\
\qquad\qquad\qquad
+\|\haep\|^2_{(L^2(\Sigma^a))^{3}} + \|\hbepp\|^2_{(L^2(\omb))^{3}} +
 \|\hbepm\|^2_{(L^2(\omb))^{3}}
=1.
\end{array}
\]
\end{rem}

%\sub
\subsection{The derivation of the scaling} 
Let us consider the possible scalings for the solution 
$\overline{U}^\ep$ and test function $U$. If, instead of a multidomain, 
one considers a single thin cylinder, the natural scaling is
(see \cite{Le2}, \cite{MuSi1}, \cite{MuSi2}):
\[
u  _\alpha (x) =    \rep\, U_\alpha 
(\rep x',x_3),  \; \; 
  u_3 (x) =    U_3
(\rep x',x_3),  \mbox{ for }x\in \Oma,    
\]
and for a single plate, it is
(see \cite{Ci1}, \cite{Le2})
\[
 u _\alpha (x) =    U_\alpha 
(  x', \ep x_3), \;  \;  
u _3 (x) =   \ep  \, U_3
(  x', \ep x_3), \mbox{ for }  x\in \Omb. 
\]
For the multidomain made of the union of the beam and the plate, 
the idea is to consider different coefficients of normalization, $\laep$ and 
$\lbep$, for $\Oaep$ and $\Obep$ respectively, that is we set 
\[\begin{array}{lll}
u^a _\alpha (x) = \laep   \rep \, U_\alpha 
(\rep x',x_3),  
& u^a_3 (x) =     \laep\, U_3
(\rep x',x_3), & \mbox{ for }x\in \Oma, \\   
u^b_\alpha (x) = \lbep   \, U_\alpha 
(  x', \ep x_3),  
& u^b _3 (x) =   \lbep \ep \, U_3
(  x', \ep x_3), &\mbox{ for }  x\in \Omb. \end{array}
\]
Then one has, with $\eaep$, $\ebep$ defined in (\ref{oeep}),
\[
e(U)(\rep x',x_3)=\frac{1}{\laep}\eaep(\uaep)(x) \mbox{ for }x \in \Oma 
\mbox{ and }
e(U)(x', \ep x_3)=\frac{1}{\lbep}\ebep(\ubep)(x) \mbox{ for }x \in \Omb 
\]
and it is easy to check that the variational equality in (\ref{varpetit}) reads, 
once each integral is written on the corresponding fixed domain,
\begin{equation}\label{eq}
\begin{array}{l}
\frac{(\rep)^2}{(\laep)^2}
\int_{\Oma}[A^a \eaep(\ouaep),\eaep(u^a)]\,dx+k^\ep \frac{\ep}{(\lbep)^2}
\int_{\Omb}[A^b \ebep(\oubep),\ebep(u^b)]\,dx = \\ \\
=\frac{1}{\laep} 
\left( \sum_{\alpha =1}^{2}
\int_{\Oma} \rep F^\ep_\alpha(\rep x',x_3)u^a_\alpha (x)\,dx +
\int_{\Oma} (\rep)^2 F^\ep_ 3(\rep x',x_3)u^a_ 3 (x)\,dx
\right) +\\\\
+\frac{1}{\lbep} 
\left( \sum_{\alpha =1}^{2}
\int_{\Omb} \ep F^\ep_\alpha( x',\ep x_3)u^b_\alpha (x)\,dx +
\int_{\Omb}  F^\ep_ 3(x',\ep x_3 )u^b_ 3 (x)\,dx
\right) +\\\\
+\frac{(\rep)^2}{\laep}\int_{\Oma}[G^\ep(\rep x',x_3),\eaep(u^a)]\,dx +
\frac{\ep}{\lbep}\int_{\Omb}[G^\ep(x',\ep x_3 ),\ebep(u^b)]\,dx +\\\\
+\frac{1}{\laep} 
\left( \sum_{\alpha =1}^{2}
\int_{\Sigma^a}  H^\ep_\alpha(\rep x',x_3)u^a_\alpha (x)\,d\sigma +
\int_{\Sigma^a } \rep  H^\ep_ 3(\rep x',x_3)u^a_ 3 (x)\,d\sigma 
\right) +\\\\
+\frac{1}{\lbep} 
\left( \sum_{\alpha =1}^{2}
\int_{\omb \setminus \repo}  H^\ep_\alpha(x',0)u^b_\alpha (x',0)\,dx'+
\int_{\omb \setminus \repo} \frac{1}{\ep} H^\ep_ 3(x',0)u^b_ 3 (x',0)\,dx'
\right) +\\\\
+\frac{1}{\lbep} \left( 
\sum_{\alpha =1}^{2}
\int_{\omb}  H^\ep_\alpha(x',-\ep)u^b_\alpha (x',-1)\,dx'+
\int_{\omb} \frac{1}{\ep} H^\ep_ 3(x',-\ep)u^b_ 3 (x',-1)\,dx' \right).  
\end{array}
\end{equation}
%As we said in the begining of Section 2, 
We decide to choose  
$\laep = \ep \lbep$, so that the junction condition written for $(u^a, u^b)$ 
reads: for almost every $x'$ in $\oma$,
\[
u^a_\alpha (x',0) = 
\frac{\laep}{\lbep}
\rep u^b_\alpha (\rep x',0) =\ep \rep u^b_\alpha (\rep x',0) \mbox{ and }
u^a_3(x',0) =  
\frac{\laep}{\lbep} \frac{1}{\ep}
u^b_3(\rep x',0) = u^b_3(\rep x',0) 
\]
(see also the begining of Section 3). Then, after dividing 
by $(\rep)^2/(\laep)^2$, writing $\lep$ instead of $\laep$, for simplicity, 
and defining the 
rescaled source terms by (\ref{fep}), (\ref{gep}), (\ref{hep}), the equality 
(\ref{eq}) is exactly the  
variational equality in (\ref{vargrand}). Finally, we recall that the 
particular choice of $\lep$ given in (\ref{lep}) makes the source terms 
bounded, but not strongly converging to zero (see also Remark 3).

\begin{rem} Since the first member of (\ref{eq}) is another way of writing 
$\int_{\Oep}[A e(\overline{U}^\ep),e(U)]\,dx$, it follows that
\begin{equation}\label{en}\begin{array}{l}
\left(\frac{\rep}{\lep}\right)^2
\left(\int_{\Oma}[A^a \eaep(\ouaep),\eaep(\ouaep )]\,dx+ \qep   
\int_{\Omb}[A^b \ebep(\oubep),\ebep(\oubep )]\,dx\right) =\\ \\
\qquad\qquad\qquad\qquad\qquad\qquad= \int_{\Oep} 
[A e(\overline{U}^\ep),e(\overline{U}^\ep)]\,dx,
\end{array}
\end{equation}
which gives the definition of the renormalized energy in (\ref{ren.en}). In 
\cite{GaMoMoMuSi}, we took $\lep =\rep$, since the initial problem 
(\ref{varpetit}) was 
supposed to be suitably normalized, {\em a priori}.
\end{rem} 

%\sub
\section{The {\em a priori} estimates and the compactness arguments}
%\noindent{\bf 
\subsection{{\em A priori} estimates:} 
%\vspace{0.3cm}

In the following, we denote by $C$ any positive constant which does not depend 
on $\ep$ and {\bf we write $\oeaep$ (resp. $\oebep$) for $e^{a \ep}(\ouaep)$ 
(resp. 
$e^{b \ep}(\oubep)$)}. Taking $u=\ouep = (\ouaep,\oubep)$ as test function in 
(\ref{vargrand}), 
we get
\begin{equation}\label{previous} \begin{array}{l}
 \int_{\Oma} [A^a \oeaep  , \oeaep  ]\, dx + \qep 
\int_{\Omb} [A^b \oebep , \oebep  ]\, dx= \\ \\
\quad \quad =\int_{\Oma}\faep . \ouaep  \, dx +\int_{\Omb}\fbep .\oubep \, dx +
 \int_{\Oma}[\gaep , \oeaep]\, dx +
\int_{\Omb}[\gbep , \oebep ]\, dx +\\ \\
\quad  \quad \quad\int_{\Sigma^a}\haep . \ouaep  \, d\sigma  + 
\int_{\omb} \left( \hbepp . \oubep_{ \;| x_3 =0}  
+ \hbepm . \oubep _{\;  | x_3 =-1} \right)\,dx'.
\end{array}
\end{equation}
From the Korn inequality, since $\ouaep$ vanishes on $T^a$ and $\oubep$ 
vanishes on $\Sigma^b$, we get for $\ep \leq 1$ and $\rep\leq 1$,
\[
\| \ouaep \|_{(H^1(\Oma))^3}  
\leq C \| e(\ouaep) \|_{(L^2(\Oma))^{3 \times 3}} 
\leq C\| \oeaep \|_{(L^2(\Oma))^{3 \times 3}}, 
\]
\[
\| \oubep \|_{(H^1(\Omb))^3}  
\leq C \| e(\oubep) \|_{(L^2(\Omb))^{3 \times 3}} 
\leq C\| \oebep \|_{(L^2(\Omb))^{3 \times 3}}, 
\]
and, by continuity of the trace mapping,
\[
\| \ouaep \|_{(L^2(\Sigma^a))^3}\leq C  \| \ouaep \|_{(H^1(\Oma))^3},
\]
\[
\| \oubep_{\; \;|x_3=0} \|_{(L^2(\omb))^3} \mbox{ and }
\| \oubep_{\; \;|x_3=-1} \|_{(L^2(\omb))^3}
\leq C  \| \ouaep \|_{(H^1(\Omb))^3}.
\]
By using the above inequalities, the coercivity of $A^a$ and $A^b$ and the 
boundedness of the source terms (see (\ref{convfa}) to (\ref{convhb}) and Remark 2), it 
follows from (\ref{previous}) that  
\[ \begin{array}{l}
C \| \oeaep \|^2_{(L^2(\Oma))^{3 \times 3}} +  C \qep 
\| \oebep \|^2_{(L^2(\Omb))^{3 \times 3}} \leq \\ \\
 \leq\left(\| \faep \|_{(L^2(\Oma))^3} + 
\| \gaep \|_{(L^2(\Oma))^{3 \times 3}} +
 \| \haep \|_{(L^2(\Sigma^a))^3}  
\right)
\| \oeaep \| _{(L^2(\Oma))^{3 \times 3}} + \\ \\
+\left(\| \fbep \|_{(L^2(\Omb))^3} +
\|\gbep \|_{(L^2(\Oma))^{3 \times 3}} +
\| \hbepp\|_{(L^2(\omb))^{3}}
+  \| \hbepm\|_{(L^2(\omb))^{3}}\right)
 \| \oebep \| _{(L^2(\Omb))^{3 \times 3}} \leq\\ \\  
\leq C \left(\| \oeaep \| _{(L^2(\Oma))^{3 \times 3}} +   
\| \oebep \| _{(L^2(\Omb))^{3 \times 3}}\right).\end{array} 
\]

$\bullet$ If $\qep$ is bounded below by some positive constant, that is if 
$q$ defined in (\ref{convq}) equals 
some positive number or $+\infty$, 
it follows that $\oeaep$ is bounded in 
$(L^2(\Oma))^{3 \times 3}$ and $\oebep$ is bounded in 
$(L^2(\Omb))^{3 \times 3}$. Then, from Korn's inequality, it results 
that $\ouaep$ is bounded in 
$(H^1(\Oma))^{3}$ and $\oubep$ is bounded in 
$(H^1(\Omb))^{3}$. Moreover, in the particular case $q=+\infty$, $\oebep$ 
tends to zero (strongly) in $(L^2(\Omb))^{3 \times 3}$ and 
$\oubep$tends to zero (strongly) in $(H^1(\Omb))^{3}$.
 
$\bullet$ Otherwise, i.e. if $\qep$ tends to zero, we define $\utep$ by
\begin{equation}
\utep = (\utaep,\utbep) = \qep \ouep = \qep(\ouaep,\oubep).
\end{equation}
It is clear that $\utep$ solves 
\begin{equation}\label{vartilde} \begin{array}{l}
\utep \in {\cal Y} ^\ep\mbox{ and } \forall u \in  {\cal Y} ^\ep, \\ \\
\quad \frac{1}{\qep}\int_{\Oma} [A^a e^{a \ep}(\utaep) , e^{a\ep}(u^a) ]\, dx +
\int_{\Omb} [A^b e^{b \ep}(\utbep) , e^{b\ep}(u^b) ]\, dx= \\ \\
\quad \quad =\int_{\Oma}\faep . u^a\, dx +\int_{\Omb}\fbep . u^b\, dx +
 \int_{\Oma}[\gaep , e^{a\ep}(u^a)]\, dx +
\int_{\Omb}[\gbep , e^{b\ep}(u^b)]\, dx +\\ \\
\quad  \quad \quad\int_{\Sigma^a}\haep . u^a\, d\sigma  + 
\int_{\omb} \left( \hbepp . u^b_{ \;| x_3 =0}  
+ \hbepm . u^b_{\;  | x_3 =-1} \right)\,dx'.
\end{array}
\end{equation}
Taking $u=\utep$ as test function in (\ref{vartilde}), it is easy to prove 
(as we have done in the case $\qep \geq C > 0$) that 
$\tilde{e}^{a\ep}:= e^{a\ep}(\utaep) =\qep \oeaep$ 
tends to zero in $(L^2(\Oma))^{3 \times 3}$, 
$\tilde{e}^{b\ep}:= e^{b\ep}(\utbep) =\qep \oebep$ is 
bounded in $(L^2(\Omb))^{3 \times 3}$, $\utaep =\qep \ouaep$ tends to zero 
in $(H^1(\Oma))^{3}$ and $\utbep =\qep \oubep$ is bounded in 
$(H^1(\Omb))^{3}$.

%\vspace{0.3cm}
%\noindent{\bf Compactness arguments:}\vspace{0.3cm}
\subsection{Compactness arguments:}
 
$\bullet$ If $\qep \rightarrow q \in (0, +\infty]$, it results from the 
a priori estimates that there exist $\ou = (\oua, \oub)$ in 
$(H^1(\Oma))^3 \times (H^1(\Omb))^3$ and $\ove =(\oea, \oeb)$ in 
$(L^2 (\Oma))^{3 \times 3}\times (L^2 (\Omb))^{3 \times 3}$, such that 
\begin{equation}\label{defou}
\ouep= (\ouaep, \oubep) \rightharpoonup \ou = (\oua, \oub)
\mbox{ weakly in }(H^1(\Oma))^3 \times (H^1(\Omb))^3,
\end{equation}
\begin{equation}\label{defoe}
\oeep= (\oeaep,\oebep) \rightharpoonup \ove= (\oea, \oeb)
\mbox{ weakly in }(L^2(\Oma))^{3 \times 3}\times (L^2(\Omb))^{3 \times 3}.
\end{equation}
Clearly $\oua =0$ on $T^a$, $\oub =0$ on $\Sigma^b$ and $\oea$, $\oeb$ are 
symmetric matrices. Moreover, from the boundedness of 
$\oeep = (\oeaep,\oebep)$ and a classical 
semicontinuity argument, we get that $\oua$ is a Bernouilli-Navier displacement and 
$\oub$ is a Kirchhoff-Love displacement:
$$
e_{\alpha \beta} (\oua) =0 \mbox{ and } e_{\alpha 3} (\oua) =0, \quad
e_{\alpha 3} (\oub) =0 \mbox{ and }e_{33} (\oub) =0,
$$
which, combined with the constraints 
$\oua =0$ on $T^a$, $\oub =0$ on $\Sigma^b$, is equivalent to (see \cite{Le1}):
\[ \begin{array}{l}
\oua \in  (H^2(0,1))^2 \times H^1(\Oma), \oua_{\alpha}(1)=
\frac{d  \oua_{\alpha}}{dx_3}(1)=0, \\ \qquad\qquad
\exists \,  \overline{\zeta}^a \in H^1(0,1), \overline{\zeta}^a (1)=0,
\oua  _3 = \overline{\zeta}^a  -
x_{1}\frac{d\, \oua_1}{d\,x_3} -
x_{2}\frac{d\, \oua_2}{d\,x_3}
 ,
\end{array} 
\]
\[
\oub \in {\cal U}^b.
%= \{u^b \in (H^1(\Omb))^2 \times H^2_0(\omb), \;
%\exists \zeta^b_\alpha \in H^1_0(\omb ), \quad 
%u^b_\alpha = \zeta^b_\alpha -
%x_3\frac{\partial u^b_3}{\partial x_\alpha}
%\}.
\]
Moreover one can prove as in \cite{MuSi2} that there exist 
$(\ova, \owa)$ 
%in ${\cal V}^a \times {\cal W}^a$ 
and $(\ovb, \owb)$ 
%in ${\cal V}^b \times {\cal W}^b$) 
such that $\oea =e^a (\oua, \ova, \owa)$ and $\oeb =e^b (\oub, \ovb, \owb)$ 
(see the definitions of $e^a$ and $e^b$ in (\ref{eab})) and such that
\[
\begin{array}{l} 
\ova \in 
	\begin{array}[t]{l}
	(H^1(\Oma))^2 \times L^2(0,1;H^1 (\oma)), \; 
	\exists \, \oc \in H^1 (0,1), \oc(1)=0, 
 	\ova_1=-\oc\,x_2, \ova_2=\oc\,x_1, 
 	\\
	\qquad \qquad\mbox{ for a.e. } 
	x_3 \in (0,1), \; 
 	\int_{\oma}\ova_3 (x',x_3)\, dx'=0, 
 	\end{array}
\\ \\
\owa  \in {\cal W}^a, \qquad 
\ovb \in {\cal V}^b, \qquad  
\owb \in {\cal W}^b 
\end{array} 
\]
and suitable expressions of $\ouaep$ (resp. $\oubep$) 
tend to $(\ova, \owa)$ (resp. $(\ovb, \owb)$). For convenience of the reader, 
the proof of this fact is given in the Appendix (see Section 8.1). In particular, $\ova$ defines 
some $\oc \in H^1(0,1)$ with $\oc(1)=0$, which is actually the limit in $L^2(0,1)$ of $\oc^\ep$ given by
\begin{equation}\label{ocep}
\oc^\ep(x_3)= \frac
{\int_{\oma}\left(x_1 \ouaep_2(x',x_3) 
-x_2 \ouaep_1(x',x_3) \right)\, dx'}
{ \rep\int_{\oma}\left(x_1^2 + x_2^2 \right)\, dx'}.
\end{equation}
To summarize, we have proved (\ref{convouep1}), (\ref{convoeep1}).

In the particular case $q=+\infty$, we have already noticed (see the 
{\em a priori} estimates) that 
\begin{equation}\label{definfini}
\oubep \rightarrow \oub = 0 \mbox{ strongly in }(H^1(\Omb))^3 \mbox{ and }
\oebep \rightarrow \oeb = 0 \mbox{ strongly in }(L^2(\Omb))^{3 \times 3},
\end{equation}
that is we have proved (\ref{convouep2}), (\ref{convoeep2}).
 
$\bullet$ If $\qep \rightarrow 0$, it results from the 
{\em a priori} estimates that 
\begin{equation}\label{def0}
%\begin{equation} 
\begin{array}{c}
\qep \, \ouaep \rightarrow 0 \mbox{ strongly in } (H^1(\Oma))^3,\;
\qep \, \oubep \rightharpoonup \oub \mbox{ weakly in } (H^1(\Omb))^3,  \\
\qep  \oeaep  \rightarrow 0 \mbox{ strongly in }(L^2(\Oma))^{3 \times 3}, \;
\qep \oebep  \rightharpoonup \oeb \mbox{ weakly in }(L^2(\Omb))^{3 
\times 3}, 
\end{array}
\end{equation}
for some $\oub \in {\cal U}^b$ and some symmetric matrix $\oeb \in 
(L^2(\Omb))^{3 \times 3}$. Again (see the Appendix, Section 8.1), there exists $(\ovb, \owb)$ in 
${\cal V}^b \times {\cal W}^b$, which are limits of suitable expressions 
of $\oubep$ and such that $\oeb=e^b(\oub, \ovb, \owb)=e^b(\ozb)$. In other words, 
we have proved (\ref{convouep3}), (\ref{convoeep3}).

%\sub
\section{The limit constraints that are due to the junction}
As for the limit constraints, it remains to prove that 

1) $\oua_\alpha (0)=0$, 

2) $\oua_ 3 (x',0) \equiv \oub_ 3 (0)$, which is equivalent to 
$\ozta(0)=\oub_ 3 (0)$ and $\frac{d \oua_\alpha}{d x_3}(0)=0$,

3) $\oc(0)=0$,\\
since the above three conditions give $\oua_{\alpha} \in (H^2_0(0,1))^2$ and 
$\oc \in H^1_0(0,1)$, so that $\oua \in {\cal U}^a$ and $\ova \in {\cal V}^a$. 
These limit constraints are derived below. 

%\sub
\subsection{Proof of $\oua_\alpha (0)=0$}

The fact that $\oua_\alpha (0)=0$ results from the following easy Lemma. 
%applied to $\mu^\ep=\ep$, $(u^{a\ep},u^{b\ep})=
%(\ouaep_{\alpha}(.,0),\oubep_{\alpha}(.,0))$ if $q\neq 0$, or to 
%$(u^{a\ep},u^{b\ep})=
%\qep(\ouaep_{\alpha}(.,0),\oubep_{\alpha}(.,0))$ if $q=0$. 
%As, 
%by compactness of the trace mapping, $\uaep(.,0)$ (or $\qep \uaep(.,0)$ 
%if $q=0$) tends to $\oua_\alpha(0)$ in $L^2(\oma)$, 
%we get $\oua_\alpha (0)=0$.
%\begin{lemma}
%Assume that $(u^{a\ep},u^{b\ep}) \in L^2(\oma)\times L^2(\omb)$ verifies 
%\[
%\mbox{ for a.e. } x' \in \oma, \qquad u^{a\ep}(x') =\mu^\ep \rep u^{b\ep}(\rep x')
%\]
%and assume that, when $\ep$ tends to zero, $\mu^\ep \rightarrow 0$ and 
%$u^{b\ep}$ is bounded in $L^2(\omb)$. Then $\uaep$ tends to zero in $L^2(\oma)$.
%\end{lemma}
\begin{lemma}
Assume that $\{u^{b\ep}\}_\ep$ is bounded in $L^2(\omb)$. Then 
$\{\rep u^{b\ep}(\rep .)\}_\ep$ is bounded in $L^2(\oma)$, for every $\oma$ 
such that $\rep \oma \subset \omb$, for any $\ep$. 
\end{lemma}
%Proof: By assumption,
%\[\begin{array}{ll}
%\int_{\oma} |u^{a\ep}(x')|^2 \, dx'&= \left(\mu^\ep \rep \right)^2
%\int_{\oma} |u^{b\ep}(\rep x')|^2 \, dx'\\ \\
%&=\left(\mu^\ep \right)^2\int_{\rep \oma} |u^{b\ep}( x')|^2 \, dx'\\ \\
%&\leq\left(\mu^\ep \right)^2\int_{\omb} |u^{b\ep}( x')|^2 \, dx',
%\end{array}
%\]
%which tends to zero with $\ep$. Hence $\uaep$ tends to zero in $L^2(\oma)$.
Proof: We have
\[
\int_{\oma} |\rep u^{b\ep}(\rep x')|^2 \, dx'=
\int_{\rep \oma} |u^{b\ep}( x')|^2 \, dx'\leq
 \int_{\omb} |u^{b\ep}( x')|^2 \, dx' \leq C.
\] 

\vspace{0.5cm}\noindent 
{\bf Application:} If $q \neq 0$, we write the junction condition for 
$\ouep_{\alpha}$:
\[  
\mbox{ for a.e. }x'\in \oma,\; \ouaep_{\alpha}(x',0)=\ep \rep \oubep_{\alpha}(\rep x',0).
\]
The first member tends to $\oua _{\alpha}(x',0)=\oua _{\alpha}(0)$ in 
$L^2(\oma)$. The second member tends to zero in this space, by Lemma 1, 
since $\oubep_{\alpha}(.,0)$ is bounded in $L^2(\omb)$, so that 
$\oua _{\alpha}(0)=0$. If $q=0$, 
the same proof applies to $\utep=\qep \ouep$.

%\sub
\subsection{Proof of $\oua_ 3 (x',0) \equiv \oub_ 3 (0)$}

This is a crucial part of this paper. 
It is derived from the following general Lemma.
%, applied to $(\uaep,\ubep)=
%(\ouaep_{3\;\;|x_3=0}, \oubep)$, if $q \neq 0$, or $(\uaep,\ubep)=\qep
%(\ouaep_{3\;\;|x_3=0}, \oubep)$, if $q = 0$. 
%(The boundedness of 
%$\{\ebep(\ubep)\}_\ep$ is given by the {\em a priori} estimates, 
%see Section 3.2.)
%\begin{lemma}

\begin{lemma}
Assume that $\ep$ and $\rep$ tend to zero, with $0<\ep^2 \ll \rep$. Let 
$ \ubep \in  (H^1(\Omb))^3$ be such that 
\[  
\ubep =0 \mbox{ on }\Sigma^b,
\]
\begin{equation}\label{bound} 
\{\ebep(\ubep)\}_\ep \mbox{ is bounded in }(L^2(\Omb))^{3 \times 3},
\end{equation}
with $\ebep$ defined in (\ref{oeep}). Then, up to a subsequence,
\begin{equation}\label{wconv}
\ubep\rightharpoonup u^b \mbox{ weakly in }(H^1(\Omb))^3,
\end{equation}
for some $u^b \in {\cal U}^b $ (in particular $u^b_3 %= u^b_3 (x') 
\in H^2_0(\omb)$). Moreover $\ubep _3(\rep .,0)$ tends to 
$u^b_ 3 (0)$ strongly in $L^2(\oma)$, for every $\oma$ 
such that $\rep \oma \subset \omb$, for any $\ep$. 
\end{lemma}
 
Proof: The first part of the lemma is classical (see \cite{Ci1}). Let us prove 
the convergence of $\ubep _3(\rep .,0)$. We define $\Uep : \omb \rightarrow \R$ by
\begin{equation}\begin{array}{ll}
\Uep (x') &=k \int_{-1}^{0}\int_{-1}^{0}\int_{t<x_3<t'}
\ubep_3(x',x_3)\,dx_3\,dt\,dt'\\ \\
&=k \int_{-1}^{0}\int_{-1}^{0} \rho (t,t') \int_{t}^{t'}
\ubep_3(x',x_3)\,dx_3\,dt\,dt',
\end{array}\end{equation} 
with $\rho (t,t')=1$ if $t<t'$, $0$ otherwise, and with $k$ chosen so that 
\begin{equation}\label{k}
k\int_{-1}^{0}\int_{-1}^{0}\rho (t,t')(t'-t)\,dt\,dt'=1.
\end{equation} 
Moreover we define $\eepa: \Omb\rightarrow \R$, $\Eepa$ and 
$\Depa:\omb \rightarrow \R$ by
\begin{equation}\label{confusion} 
\eepa=2 \, e_{\alpha 3} (\ubep)=
\frac{\partial\ubep_\alpha}{\partial x_3}+
\frac{\partial\ubep_3}{\partial x_\alpha},
\end{equation} 
\begin{equation}\label{confusionbis} 
\Eepa (x')=k \int_{-1}^{0}\int_{-1}^{0} \rho (t,t') \int_{t}^{t'}
\eepa(x',x_3)\,dx_3\,dt\,dt',
\end{equation}
\begin{equation}\label{defD}\begin{array}{ll} 
\Depa (x')&=k \int_{-1}^{0}\int_{-1}^{0} \rho (t,t') \int_{t}^{t'}
\frac{\partial\ubep_\alpha}{\partial x_3} (x',x_3)\,dx_3\,dt\,dt'\\ \\
&=k \int_{-1}^{0}\int_{-1}^{0} \rho (t,t')
\left(
\ubep_\alpha (x',t')-\ubep_\alpha (x',t)\right)\,dt\,dt'.
\end{array}\end{equation}
It is clear that
\begin{equation}\label{1} 
\nabla \Uep = \Eep - \Dep.
\end{equation}
Still denoting by $C$ various constants that do not depend on $\ep$, we have 
from Cauchy-Schwarz inequality:
$$ 
|\Eepa (x')|^2 
\leq C\int_{-1}^{0}|\eepa(x',x_3) |^2\,dx_3,
$$ 
which gives, by definition of $\eepa$ and by (\ref{bound}),
\begin{equation}\label{2} 
\|\Eepa\| _{L^2(\omb)}\leq C \|\eepa\| _{L^2(\Omb)} = 
C \| e_{\alpha 3}(\ubep)\|_{L^2(\Omb)}  \leq C\ep.
\end{equation}
From (\ref{defD}), Cauchy-Schwarz inequality 
and the boundedness of $\ubep_\alpha$ in $H^1_0(\Omb)$, we have
\begin{equation}\label{3}
\|\Dep_{\alpha}\| _{H^1_0(\omb)}\leq C \|\ubep_{\alpha}\| _{H^1_0(\Omb)}\leq C. 
\end{equation}
From (\ref{1}), we get the following decomposition:
$$
\Uep = {\hat U}^\ep+{\tilde U}^\ep,
$$
with ${\hat U}^\ep$, ${\tilde U}^\ep$ the respective solutions 
in $H^1_0(\omb)$ of 
$$
-\Delta {\hat U}^\ep = -div \Eep \quad  \mbox{ and }
\quad -\Delta {\tilde U}^\ep = 
div \Dep\quad \mbox{ in }\omb,
$$
and from (\ref{2}), (\ref{3}), 
\begin{equation}\label{oubli}
\|\nabla {\hat U}^\ep\| _{(L^2(\omb))^2}\leq \|\Eep\|_{(L^2(\omb))^2} \leq C \ep,
\end{equation}
\begin{equation}\label{num1}
{\hat U}^\ep \rightarrow 0 \mbox{ in }H^1_0(\omb),
\end{equation}
\begin{equation}\label{num2} 
\| {\tilde U}^\ep\| _{H^2(\omb)}\leq C \| div \Dep\|_{L^2(\omb)} \leq C.
\end{equation}  
But, using (\ref{wconv}) and (\ref{k}), it is easy to prove that
$$ 
\Uep\rightharpoonup u^b_3 = u^b_3(x') \mbox{ weakly in }L^2(\omb),
$$ 
which gives, by virtue of (\ref{num1}), (\ref{num2}),
$$
{\tilde U}^\ep = U^\ep - {\hat U}^\ep\rightharpoonup  u^b_3 \mbox{ weakly in }H^2(\omb).
$$
Then, as the embedding $H^2(\omb) \subset {\cal C}^{0}(\overline\omb)$ 
is compact,
for $\omb$ bidimensional, we get that 
\begin{equation}\label{unif}
{\tilde U}^\ep \rightarrow u^b_3 \mbox{ in }{\cal C}^{0}(\overline\omb).
\end{equation}
This is enough to prove that $\ubep_3 (\rep .,0)$ tends to $u^b_3(0)$ 
strongly in $L^2(\oma)$. 

Actually we have, for a.e. $x'$ in $\oma$,
\begin{equation}\label{brackets}\begin{array}{ll} 
\ubep_3 (\rep x',0) - u^b_ 3 (0) 
&=\left[\ubep_ 3 (\rep x',0) - \Uep(\rep x')\right]+\\
&+\left[\Uep(\rep x') - {\tilde U}^\ep (\rep x')\right]+\\
&+\left[{\tilde U}^\ep (\rep x') -u^b_3(\rep x')  \right]+\\
&+\left[u^b_3(\rep x')  - u^b_ 3 (0)\right].
\end{array}
\end{equation}
We are going to show that each of the above brackets tends to 
zero, strongly in $L^2(\oma)$.

As for the first bracket, we have:
\begin{equation}\label{third}
\int_{\oma} |\ubep_ 3 (\rep x',0) - \Uep(\rep x')|^2 \, dx'=\frac{1}{(\rep)^2}
\int_{\rep \oma}|\ubep_ 3 (x',0) - \Uep(x')|^2 \, dx'.
\end{equation}
But, by using (\ref{k}),
\[\begin{array}{ll} 
\Uep(x') &= k \int_{-1}^{0}\int_{-1}^{0}\rho(t,t')\int_{t}^{t'}
\left(
\ubep_3(x',0)+\int_{0}^{x_3} \frac{\partial \ubep_3}{\partial x_3}(x',y_3)\, dy_3
\right)\, dx_3\,dt\,dt'=\\
&=\ubep_3(x',0) + k \int_{-1}^{0}\int_{-1}^{0}\rho(t,t')\int_{t}^{t'}
\int_{0}^{x_3} \frac{\partial \ubep_3}{\partial x_3}(x',y_3)\, dy_3\, 
dx_3\,dt\,dt',
\end{array}
\]
so that
\[
|\Uep(x')-\ubep_3(x',0)| \leq C \int_{-1}^{0} 
\left|\frac{\partial \ubep_3}{\partial x_3}(x',y_3)\right| \,dy_3 ,
\]
and this gives with (\ref{bound})
\[
\int_{\rep \oma}|\ubep_3(x',0)-\Uep(x')|^2 \, dx' \leq C \int_{\Omb}
\left|\frac{\partial \ubep_3}{\partial x_3}\right|^2 \, dx \leq C \ep^4.
\]
Coming back to (\ref{third}), it results that  
\[  
\int_{\oma} |\ubep_ 3 (\rep x',0) - \Uep(\rep x')|^2 \, dx'\leq 
C \frac{\ep^4}{(\rep)^2},
\]
which tends to zero, since we have assumed that $\ep^2 \ll \rep$. 
Now we consider the second bracket in (\ref{brackets}), that is 
${\hat U}^\ep(\rep x')$, and we are going to prove that its $L^2-$norm 
tends to zero, again if $\ep^2 \ll \rep$. In fact, from 
Cauchy-Schwarz inequality, the continuity of the imbedding 
$H^1_0(\omb) \subset L^4(\omb)$ (actually $L^q(\omb)$, for every finite $q$, 
in dimension 2) and from (\ref{oubli}),
\[\begin{array}{ll}  
\int_{\oma}|{\hat U}^\ep(\rep x')|^2\, dx'
&=\frac{1}{(\rep)^2}\int_{\rep\oma}|{\hat U}^\ep(  x')|^2\, dx'\\
&\leq\frac{1}{(\rep)^2} 
\left(\int_{\rep\oma}|{\hat U}^\ep(  x')|^4\, dx'   \right)^{\frac{1}{2}}
|\rep\oma|^{\frac{1}{2}}\\
&\leq C(\rep)^{-1} \| {\hat U}^\ep \|_{L^4(\omb)}^2\\
%&\leq \mbox{ (by continuity of the imbedding }H^1_0(\omb) \subset L^4(\omb))\\
&\leq C (\rep)^{-1} \| {\hat U}^\ep \|_{H^1_0(\omb) }^2\\
%&\leq \mbox{ (see (\ref{oubli}))} \\
&\leq C \frac{(\ep)^2}{\rep}.
\end{array}
\]
The third bracket in (\ref{brackets}) tends to zero with $\ep$, 
in $L^\infty-$norm, and also the fourth one tends to zero, by virtue 
of (\ref{unif}). 
This concludes the proof of Lemma 2.

\vspace{0.5cm}\noindent 
{\bf Application:} If $q \neq 0$, we write the junction condition for $\ouep_3$:
\[  
\mbox{ for a.e. }x'\in \oma,\;\ouaep_{3}(x',0)= \oubep_{3}(\rep x',0).
\]
The first member tends to $\oua _{3}(x',0)$ in 
$L^2(\oma)$. The second member tends to $\oub _{3}(0)$ in this space, 
by Lemma 2. It follows that, for a.e. $x'$ in $\oma$, 
$\oua _{3}(x',0)=\oub _{3}(0)$. If $q=0$, 
the same proof applies to $\utep=\qep \ouep$.

%\sub
\subsection{Proof of $\oc(0)=0$}
This part also is crucial.  

\begin{lemma}
Assume that $\ep$ and $\rep$ tend to zero, with $0<\ep^2 \ll \rep$. Let 
$(\uaep,\ubep) \in (H^1(\Oma))^{3}\times (H^1(\Omb))^{3}$ 
be such that
\begin{equation}
\uaep_{\;\;|x_3=1}=0,
\end{equation}
\begin{equation}\label{bc}
\mbox{a.e. }x' \in \oma, \; \uaep_{\alpha} (x',0)= 
\ep \rep \ubep_{\alpha} (\rep x',0) ,
\end{equation}
\begin{equation}\label{bound2} 
\{\eaep(\uaep) \}_\ep \mbox{ is bounded in }(L^2(\Oma))^{3 \times 3},
\end{equation}
\begin{equation}\label{bound3}
\{\ubep_{\alpha}\}_\ep \mbox{ is bounded in }H^1(\Omb). 
\end{equation}
Let $c^\ep$ be defined by 
\begin{equation}\label{cep}
c^\ep(x_3)= \frac
{\int_{\oma}\left(x_1 \uaep_2(x',x_3) 
-x_2 \uaep_1(x',x_3) \right)\, dx'}
{ \rep\int_{\oma}\left(x_1^2 + x_2^2 \right)\, dx'}.
\end{equation}
Then $c^\ep \rightarrow c$ in $L^2(0,1)$, for some $c$ in $H^1_0(0,1)$.
\end{lemma}

Proof: For $\alpha=1,2$, we define $x^R_\alpha$ by $x^R_1=-x_2$, $x^R_2=x_1$ 
and we set
\[
\vaep=\frac{\uaep}{\rep},
\]
\[
\eepa=2e_{\alpha 3}(\vaep)=2e^{a\ep}_{\alpha 3}(\uaep)
\]
(without confusion with $\eepa$ appearing in (\ref{confusion})),
\[
\mepa=\frac{1}{|\oma|}\int_{\oma}\vaep_\alpha\, dx',
\]
\[
\rho^\ep_\alpha=\frac{1}{\rep}\left[\vaep_\alpha -c^\ep x^R_\alpha -
\mepa \right],
\]
with $c^\ep$ given by (\ref{cep}).

We begin by giving {\bf two {\em a priori} estimates}. Due to (\ref{bound2}), we have
\begin{equation}\label{bounde}
\left\|e^\ep\right\|_{(L^2(\Oma))^2} \leq C .
\end{equation}
As for $\rho^\ep$, it follows from (\ref{conddom}) that $\rho^\ep_\alpha(.,x_3)$ 
has mean value zero on $\oma$, for every $x_3$ and, as 
$e_{\alpha \beta}(\rho^\ep)=(1/\rep)e_{\alpha \beta}(\vaep)$, we get from 
the Poincar\'e-Wirtinger inequality for elasticity
\[\begin{array}{ll}
\left\|\rho^\ep\right\|^2_{(L^2(\Oma))^2} &\leq C \sum_{\alpha\beta}
\left\|e_{\alpha \beta}(\rho ^\ep) \right\|^2_{L^2(\Oma)}
=\frac{C}{(\rep)^2}\sum_{\alpha\beta}
\left\|e_{\alpha \beta}(\vaep) \right\|^2_{L^2(\Oma)}\\ \\
&= C  \sum_{\alpha\beta}
\left\|\eaep_{\alpha \beta}(\uaep) \right\|^2_{L^2(\Oma)},
\end{array}
\]
which gives, with (\ref{bound2}),
\begin{equation}\label{boundrho}
\left\|\rho ^\ep\right\|_{(L^2(\Oma))^2} \leq C .
\end{equation}

Now we prove that {\bf one can derive a single equation}, of the form 
$c^\ep =K^\ep-\rep R^\ep$, {\bf from the system of two equations 
$c^\ep x^R_\alpha + \mepa=\vaep_\alpha -\rep \rho^\ep_\alpha$}. 
(This is a very tricky argument appearing in \cite{MuSi2}, 
see also Section 8.1.) Indeed we get by differentiating the previous system
with respect to $x_3$,  
\begin{equation}\label{..} 
\frac{dc^\ep}{dx_3}x^R_\alpha + \frac{d \mepa}{dx_3} + 
\frac{\partial \vaep_3}{\partial x_\alpha}=\eepa -\rep 
\frac{\partial\rho^\ep_\alpha}{\partial x_3}, \quad \forall \alpha=1,2 .
\end{equation}
After multiplying (\ref{..}) by a test function $\varphi_\alpha \in {\cal D}(\oma)$, 
summing over $\alpha$ and integrating over $\oma$, we have
\begin{equation}\label{...}\begin{array}{l}
\frac{dc^\ep}{dx_3}\int_{\oma}\sum_{\alpha}\varphi_\alpha  x^R_\alpha \, dx'
+\sum_{\alpha}\frac{d \mepa}{dx_3}\int_{\oma}\varphi_\alpha   \, dx'-
\int_{\oma}\vaep_3 div \varphi\, dx'=\\ \\ \qquad \qquad\qquad \qquad=
\int_{\oma}\sum_{\alpha}\eepa \varphi_\alpha\, dx' - \rep \frac{d}{dx_3}
\int_{\oma}\sum_{\alpha}\rho^\ep_{\alpha}\varphi_\alpha\, dx'.
\end{array}
\end{equation}
We choose the test function $\varphi_\alpha$ so that 
\begin{equation}\label{test1}
\int_{\oma}\sum_{\alpha}\varphi_\alpha  x^R_\alpha \, dx' =1
\end{equation}
\begin{equation}\label{test2}
\int_{\oma}\varphi_\alpha   \, dx'=0, \quad\forall \alpha =1,2,
\end{equation}
\begin{equation}\label{test3}
div \varphi =0.
\end{equation}
It is easy to check that such test function does exist: take e.g. 
$$\varphi_1= \frac{\partial \phi}{\partial x_2}, 
\varphi_2= -\frac{\partial \phi}{\partial x_1}, \mbox{ with }
\phi \in {\cal D}(\oma), \int_{\oma}\phi \, dx' = \frac{1}{2}.$$
Now we set (without confusion with $E^\ep$ appearing in (\ref{confusionbis}))
\[ 
E^\ep=\int_{\oma}e^\ep.\varphi\, dx', 
\qquad K^\ep = - \int_{x_3}^{1}E^\ep(y_3)\, dy_3,
\qquad R^\ep=\int_{\oma}\rho^\ep.\varphi\, dx',
\]
where "." denotes the scalar product in $\R^2$. Then (\ref{...}) reads:
\[ 
\frac{dc^\ep}{dx_3}=\frac{d K^\ep}{dx_3}  -\rep \frac{d R^\ep}{dx_3}, 
\]
which gives by integration
\begin{equation}\label{cep'}
c^\ep =K^\ep-\rep R^\ep,
\end{equation}
since $c^\ep(1)=K^\ep(1)=0$ and since also $R^\ep(1)=0$, because 
$\rho^\ep(1)=0$.

Now {\bf we pass to the limit} in (\ref{cep'}). Due to (\ref{bounde}) 
and (\ref{boundrho}), $E^\ep$ and $R^\ep$ 
are bounded in $L^2(0,1)$. Moreover it follows 
that $K^\ep$ is bounded in $H^1(0,1)$, 
since by Poincar\'e inequality, 
\[ 
\left\|K^\ep\right\|^2_{H^1(0,1)} \leq 
C \displaystyle\int_0^1 \left| \frac{dK^\ep}{dx_3} \right|^2\, dx_3 
=C \displaystyle\int_0^1 \left|  E^\ep \right|^2\, dx_3  \leq C.
\]   
Then we deduce that there exists $c$ in $H^1(0,1)$, with $c(1)=0$, such that 
\[
K^\ep \rightharpoonup c \mbox{ weakly in }H^1(0,1), \mbox{ hence } 
K^\ep \rightarrow c \mbox{ strongly in }{\cal C}^0(0,1).
\]
As $\rep R^\ep$ tends to zero strongly in $L^2(0,1)$, it follows 
from (\ref{cep'}) and the above that $c^\ep$ 
tends to $c$ strongly in $L^2(0,1)$. 

It remains to {\bf prove that $c$ vanishes at the origin}. For this, we notice 
that $K^\ep(0)\rightarrow c(0)$. But one has also another expression for 
$K^\ep$. Actually, from (\ref{cep'}), (\ref{test1}), (\ref{test2}) 
and by definition of $\rho^\ep_{\alpha}$, 
\[\begin{array}{ll}
K^\ep=\rep R^\ep +c^\ep&=\rep\int_{\oma}\rho^\ep.\varphi\, dx'+c^\ep
\int_{\oma}\sum_{\alpha}\varphi_\alpha  x^R_\alpha \, dx'+
\sum_{\alpha}\mepa\int_{\oma}\sum_{\alpha}\varphi_\alpha  \, dx'\\ \\
&=\sum_{\alpha}\int_{\oma}\left(\rep\rho^\ep_{\alpha}+
c^\ep x^R_\alpha +\mepa
  \right)\varphi_\alpha\, dx'=
\sum_{\alpha}\int_{\oma}\vaep_{\alpha}\varphi_\alpha \, dx',
\end{array}
\]
that is 
\[
K^\ep=\sum_{\alpha}\int_{\oma}\vaep_{\alpha}\varphi_\alpha \, dx'
=\sum_{\alpha}\int_{\oma}\frac{1}{\rep}\uaep_{\alpha}\varphi_\alpha \, dx'
\]
and in particular, due to the boundary condition (\ref{bc}),
\[
K^\ep(0)=\ep\sum_{\alpha}\int_{\oma} \ubep_{\alpha}(\rep x',0)\varphi_\alpha
 \, dx'.  
\]
Hence, by using H\"{o}lder inequality, the continuity of the imbedding of 
$H^{\frac{1}{2}}(\omb)$ in $L^4(\omb)$ (in dimension 2), 
the continuity of the trace mapping from $H^1(\Omb)$ to $H^{\frac{1}{2}}(\omb)$ 
and (\ref{bound3}), we get
\[\begin{array}{ll}
|K^\ep(0)| &\leq C \ep\sum_{\alpha}\int_{\oma} |\ubep_{\alpha}(\rep x',0)|
 \, dx' = C \frac{\ep}{(\rep)^2}
\sum_{\alpha}\int_{\rep\oma} |\ubep_{\alpha}(x',0)| \, dx'\\ \\
&\leq C \frac{\ep}{(\rep)^2}\sum_{\alpha}
\left[\int_{\rep\oma} |\ubep_{\alpha}(x',0)|^4 \, dx'\right]^{\frac{1}{4}}
|\rep\oma|^{\frac{3}{4}} \\ \\
&=C\ep(\rep)^{-\frac{1}{2}}\sum_{\alpha}
\left[\int_{\rep\oma} |\ubep_{\alpha}(x',0)|^4 \, dx'\right]^{\frac{1}{4}}\\ \\
&\leq C\ep(\rep)^{-\frac{1}{2}}\sum_{\alpha}
\left\| \ubep_{\alpha}(.,0) \right\|_{L^4(\omb)}
\leq C\ep(\rep)^{-\frac{1}{2}}\sum_{\alpha}
\left\| \ubep_{\alpha}(.,0) \right\|_{H^{\frac{1}{2}}(\Omb)} 
\leq C\ep(\rep)^{-\frac{1}{2}},
\end{array} 
\]
which tends to zero, since $0<\ep^2 \ll \rep$. As we have proved that 
$K^\ep(0)\rightarrow c(0)$, we conclude that $c(0)=0$.

%\sub
\section{The use of convenient test functions}

This is the third crucial part of the paper, at least in the case 
$q\in (0,+\infty)$.

%\sub
\subsection{The case $q=+\infty$}
Remark that ${\cal Z}_\infty = \{z^a \in {\cal Z}^a, 
\zeta^a(:=\zeta^a(u^a)) \in H^1_0(0,1) \}$. Let $u^a \in {\cal U}^a$, with 
$\zeta^a\in H^1_0(0,1)$ and let $(v^a, w^a)$ be such that
\begin{equation}\label{reg1}\begin{array}{l} 
v^a_1 =-cx_2 \mbox{ and }v^a_2 =cx_1, \mbox{ with }c\in H^1_0(0,1);\quad 
v^a_3 \in {\cal C}^1(\overline{\Oma}), v^a_{3\;|x_3=0}=0;\\ \\
w^a_{\alpha} \in {\cal C}^1(\overline{\Oma}), w^a_{\alpha \;|x_3=0}=0; 
\quad w^a_3 =0.
\end{array}
\end{equation}
In other words, $v^a$ and $w^a$ satisfy all the conditions given in the 
definitions of ${\cal V}^a$ and ${\cal W}^a$, but the integral ones; 
moreover $v^a_3$ and $w^a_\alpha$ belong to 
\begin{equation}\label{r1}
{\cal R} =\{v \in {\cal C}^1(\overline{\Oma}),\,v_{|x_3=0} = 0  \}.
\end{equation}
Let $\uaep = u^a + \rep v^a + (\rep)^2 w^a$. Then it is easy to check that 
$u=(\uaep,0)$ is in ${\cal Y}^\ep$. Taking it 
as test function in the variational equation 
of Problem (\ref{varpetit}), we get
\begin{equation}\label{341}  
\int_{\Oma} [A^a  \oeaep , e^{a\ep}(\uaep) ]\, dx = 
\int_{\Oma}\faep . \uaep \, dx  +\int_{\Oma}[\gaep , e^{a\ep}(\uaep )]\, dx +
 \int_{\Sigma^a}\haep . \uaep \, d\sigma.
\end{equation}
But we have, since 
$e_{\alpha \beta}(u^a)=e_{\alpha 3}(u^a)= e_{\alpha \beta}(v^a)
=e_{33}(w^a)=0$,
\[
 e^{a\ep}(\uaep) = \left(
\begin{array}{cc}  e_{\alpha \beta}(w^a) 
&  e_{\alpha  3}(v^a ) \\ \\
 e_{3 \alpha}( v^a  )
& e_{33}(u^a )\end{array}\right) + \rep
\left(
\begin{array}{cc}   0
&  e_{\alpha  3}(w^a ) \\ \\
 e_{3 \alpha}( w^a  )
& e_{33}(v^a )\end{array}\right),
\]
so that $e^{a\ep}(\uaep)$ tends to $e^a(z^a)=e^a(u^a,v^a,w^a)$ (strongly) in 
$(L^2(\Oma))^{3 \times 3}$. 
Moreover $\uaep$ tends to $u^a$ (strongly) in $(H^1(\Oma))^{3}$ and we have 
seen in (\ref{defoe}) that 
$ \oeaep$ tends to $e^a(\oza)$ 
weakly in $(L^2(\Oma))^{3 \times 3}$. By passing to the limit in (\ref{341}), 
using (\ref{convfa}), (\ref{convga}), (\ref{convha}), it follows that 
\begin{equation}\label{ve}
\int_{\Oma} [A^a e^{a}(\oza) , e^{a}(z^a ) ]\, dx = 
\int_{\Oma}f^a . u^a \, dx  +\int_{\Oma}[g^a  , e^{a}(z^a )]\, dx +
 \int_{\Sigma^a} h^a . u^a   \, d\sigma,
\end{equation}
which is the variational equation of (\ref{vargrand2}). It holds also true, 
by density and continuity,  for every $(v^a, w^a)$ 
such that
\[\label{reg1}\begin{array}{l} 
v^a_1 =-cx_2 \mbox{ and }v^a_2 =cx_1, \mbox{ with }c\in H^1_0(0,1);\quad 
v^a_3 \in L^2(0,1;H^1(\oma));\\ \\
w^a_{\alpha} \in L^2(0,1;H^1(\oma)); 
\quad w^a_3 =0,
\end{array}
\]
i.e. for every $(v^a, w^a)$ satisfying the conditions given in the 
definitions of ${\cal V}^a$ and ${\cal W}^a$, but the integral ones. 
(Note that ${\cal R}$ defined in (\ref{r1}) is dense in 
$L^2(0,1;H^1(\oma))$.) In particular (\ref{ve}) is also true for any 
$z^a\in {\cal Z}_\infty$. This means that $\oza$ solves the variational 
problem (\ref{vargrand2}).

%\sub
\subsection{The case $q=0$}
Then we have seen that $\utep = \qep \ouep = \qep(\ouaep,\oubep)$ 
solves (\ref{vartilde}) and that (see (\ref{def0})) $\utaep =\qep \ouaep$ tends to $\oua =0$
in $(H^1(\Oma))^{3}$, $\tilde{e}^{b\ep} =\qep \oebep$ 
tends to $\oeb=e^b(\ozb)$ weakly in $(L^2(\Omb))^{3 \times 3}$, for some 
$\ozb$ in ${\cal Z}_0$. (In particuler $\oub_3(0)=0$.) 

Let $B$ be some given small ball, with center zero, in $\omb$. 
Let $z^b=(u^b, v^b, w^b)$ be such that 
\begin{equation}\label{reg2}
\begin{array}{l}
u^b \in {\cal U}^b, \quad\zeta^b_\alpha(:=\zeta^b_\alpha(u^b))
\equiv u^b_3\equiv 0 \mbox { in } B; \\ \\
v^b_\alpha \in {\cal C}^1 (\overline{\Omb}), 
v^b_{\alpha}\equiv 0 \mbox { in } B \times \{0\};\quad v^b_3 =0; \\ \\
w^b_{\alpha}=0; \quad w^b_{3}\in {\cal C}^1 (\overline{\Omb}),\;
w^b_{3 }\equiv 0 \mbox { in } B\times \{0\}.
\end{array}
\end{equation}
In other words, $z^b$ satisfies all the conditions given in the 
definition of ${\cal Z}_0$, 
except the integral ones; 
moreover $\zeta^b_\alpha$ and $u^b_3$ vanish in $B$, 
$v^b_\alpha$ and $w^b_3$ belong to 
${\cal C}^1(\overline{\Omb})$ and vanish in $B \times\{ 0\}$.
Let $\ubep = u^b + \ep v^b + (\ep)^2 w^b$. Then it is easy to check that 
$u=(0, \ubep)$ is in ${\cal Y}^\ep$, for $\ep$ small enough.
Taking it 
as test function in the variational equation 
of Problem (\ref{varpetit}), we get
\begin{equation}\label{342}\begin{array}{ll}  
\int_{\Omb} [A^b \tilde{e}^{b\ep}, e^{b\ep}(\ubep )]\, dx &= 
\int_{\Omb}\fbep . \ubep \, dx  +\int_{\Omb}[\gbep , e^{b\ep}(\ubep )]\, dx \\ \\
&+
\int_{\omb} \left( \hbepp . \ubep_{ \;| x_3 =0}  
+ \hbepm . \ubep _{\;  | x_3 =-1} \right)\,dx'.\end{array} 
\end{equation}
But we have, since 
$e_{\alpha 3}(u^b)=e_{33}(u^b)= e_{33}(v^b)
=e_{\alpha \beta }(w^b)=0$,
\[
 e^{b\ep}(\ubep) = \left(
\begin{array}{cc}  e_{\alpha \beta}(u^b) 
&  e_{\alpha  3}(v^b ) \\ \\
 e_{3 \alpha}( v^b )
& e_{33}(w^b )\end{array}\right) + \ep
\left(
\begin{array}{cc} e_{\alpha \beta}(v^b)  
&  e_{\alpha  3}(w^b ) \\ \\
 e_{3 \alpha}( w^b  )
&  0\end{array}\right),
\]
so that $e^{b\ep}(\ubep)$ tends to $e^b(z^b)$ (strongly) in 
$(L^2(\Omb))^{3 \times 3}$. 
Moreover $\ubep$ tends to $u^b$ (strongly) in $(H^1(\Omb))^{3}$ and we have 
seen that 
$\tilde{e}^{b\ep}$ tends to $e^b(\ozb)$ 
weakly in $(L^2(\Omb))^{3 \times 3}$. By passing to the limit in (\ref{342}), 
using (\ref{convfb}), (\ref{convgb}), (\ref{convhb}), it follows that 
\begin{equation}\label{ve2} \begin{array}{ll}
\int_{\Omb} [A^b e^{b}(\ozb) , e^{b}(z^b ) ]\, dx &= 
\int_{\Omb}f^b . u^b \, dx  +\int_{\Omb}[g^b  , e^{b}(z^b )]\, dx \\ \\
&+\int_{\omb} \left( h^b_+. u^b_{ \;| x_3 =0}  
+  h^b_ - . u^b  _{\;  | x_3 =-1} \right)\,dx'  , \end{array}
\end{equation}
for every $z^b=(u^b,v^b,w^b)$ having the regularity of (\ref{reg2}).

But the following density results are proved in the Appendix (Section 8.2). 
First, from Lemma 5, any $\zeta^b_{\alpha} \in H^1_0(\omb)$ can be approached 
(in $H^1_0(\omb)$-norm) by a sequence $\zeta^{bn}_{\alpha}$, with 
$\zeta^{bn}_{\alpha}\equiv 0$ in a ball $B^n$ of radius $r^n$, tending to zero. 
Moreover, from Lemma 6, any $u^b_{3} \in H^2_0(\omb)$, 
with $u^b_{3}(0)=0$ can be approached 
(in the weak topology of $H^2_0(\omb)$) by a sequence 
$u^{bn}_{3}$ that vanishes in $B^n$. 
Finally, from Lemma 7, any $v^b_{\alpha}$ (or $w^b_{3}$) in 
$L^2(\omb;H^1(-1,0))$ can be approached (in $L^2(\omb;H^1(-1,0))$-norm) 
by a sequence of functions $v^{nb}_{\alpha}$ (or $w^{nb}_{3}$) in 
${\cal C}^1 (\overline{\Omb})$ that vanish in $B^n \times 
\{0\}$. 
 
By continuity, it results that (\ref{ve2}) holds true for any 
\[
\begin{array}{l}
u^b \in {\cal U}^b;\quad u^b_3(0)=0;\\ \\
v^b_\alpha \in  L^2(\omb;H^1(-1,0));\quad v^b_3 =0; \\ \\
w^b_{\alpha}=0; \quad w^b_{3}\in L^2(\omb;H^1(-1,0)),
\end{array}
\]
i.e. for every $z^b$ satisfying the conditions given in the 
definitions of ${\cal Z}_0$, but the integral ones. 
In particular (\ref{ve2}) is also true for any 
$z^b\in {\cal Z}_0$. This means that $\ozb$ solves the variational 
problem (\ref{vargrand3}).

%\sub
\subsection{The case $q\in (0,+\infty)$}
Let $z=(z^a,z^b)=((u^a,v^a,w^a),(u^b,v^b,w^b))\in 
({\cal C}^1(\overline{\Oma}))^3 \times 
({\cal C}^1(\overline{\Omb}))^3 $. We assume that $z$ satisfies all the 
conditions required in the definition of $Z$, except the integral ones, and we assume that it is "more regular". 
In particular $u^a_3(x',0)\equiv u^b_3(0)$, that is $\zeta^a(0)=u^b_3(0)$. 
The precise requirements are given below:
\begin{equation}\label{reg3}\begin{array}{l} 
u^a \in {\cal U}^a, \;u^a_{\alpha} \in{\cal C}^2[0,1], 
\;\zeta^a\in{\cal C}^1[0,1]; \\ \\
v^a_1 = -cx_2\mbox{ and }v^a_2 = cx_1\mbox{ with } c\in{\cal C}^1[0,1], 
c(0)=c(1)=0;\; \\
v^a_3 \in{\cal C}^1(\overline{\Oma}), v^a_{3\; |x_3=0}=0;\\ \\
w^a_{\alpha} \in{\cal C}^1(\overline{\Oma}),
w^a_{\alpha\; |x_3=0}=0; w^a_3=0; \\ \\
u^b \in {\cal U}^b, \zeta^b_{\alpha}\in{\cal C}^1(\overline{\omb})
\cap H^1_0(\omb);\;  u^b_{3}\in{\cal C}^1(\overline{\omb})
\cap H^2_0(\omb);\\ \\
v^b \in ({\cal C}^1(\overline{\Omb}))^2 \times \{0\}; \; 
w^b \in (\{0\})^2 \times{\cal C}^1(\overline{\Omb}); \; 
v^b_{\alpha}\mbox{ and }w^b_3 =0 \mbox{ on }\Sigma^b;\\ \\
u^b_3(0)=\zeta^a(0).  
\end{array}
\end{equation}
 
We are going to define a convenient test function $u^\ep=(\uaep,\ubep)$ in 
${\cal Y}^\ep$. 

$\bullet$ In $\Omb$, we choose
\begin{equation}\label{ubep}
\ubep=u^b+\ep v^b+\ep^2 w^b.
\end{equation}
As the couple $u^\ep=(\uaep,\ubep)$ needs to satisfy the transmission 
conditions (\ref{tc}), i.e. $$\mbox{ for a.e. }x' \in \oma, \;
u^a_\alpha (x',0) = \ep \rep u^b_\alpha (\rep x',0) \mbox{ and }\;
u^a_3(x',0) =  u^b_3(\rep x',0) \},$$ then, necessarily, $\uaep(x',0)$ 
is given by:
\[
\uaep_\alpha (x',0) = \ep \rep \left( \zeta^b _\alpha (\rep x') 
+\ep v^b _\alpha (\rep x',0) \right),
\]
\[
\uaep_3 (x',0) = u^b_3(\rep x') + \ep^2 w^b_3(\rep x',0).
\]

$\bullet$ In $\Oma \cap \{x_3>\rep\}$, we choose
\begin{equation}\label{uaep1}
\uaep=u^a+\rep v^a+(\rep)^2 w^a.
\end{equation}

$\bullet$ In $\Oma \cap \{0<x_3<\rep\}$, $\uaep$ is obtained by 
linear interpolation between $\uaep(x',0)$ and $\uaep(x',\rep)$:
\[
\uaep (x',x_3)=\frac{x_3}{\rep} \left(u^a(x',\rep) +\rep v^a(x',\rep) +
(\rep)^2 w^a(x',\rep)  \right) +\left(1-\frac{x_3}{\rep}\right)
\uaep(x',0),
\]
that is (see above)
\begin{equation}\label{uaep2}\begin{array}{ll}
\mbox{for }0<x_3<\rep, \quad \uaep_{\alpha} (x',x_3)
&=  \displaystyle \frac{x_3}{\rep} \left(u^a_{\alpha}(\rep) 
+\rep v^a_{\alpha}(x',\rep) +
(\rep)^2 w^a_{\alpha}(x',\rep)  \right)+ \\ \\
&+\left(1-\frac{x_3}{\rep}\right)
\ep \rep\left(\zeta^b_{\alpha}(\rep x')+\ep v^b_{\alpha}(\rep x',0)\right),
\end{array}
\end{equation}
\begin{equation}\label{uaep3}\begin{array}{ll}
\mbox{for }0<x_3<\rep, \quad 
\uaep_{3} (x',x_3)
&= \displaystyle   \frac{x_3}{\rep} \left(u^a_{3}(x',\rep) 
+\rep v^a_{3}(x',\rep)    \right)+ \\ \\
&+\left(1-\frac{x_3}{\rep}\right)
 \left(u^b_3(\rep x')+\ep^2 w^b_{3}(\rep x',0)\right).
\end{array}
\end{equation}
Taking $u^\ep=(\uaep,\ubep)$ as test function in the variational equation 
of Problem (\ref{vargrand}), we get
\begin{equation}\label{test} \begin{array}{l}
\int_{\Oma} [A^a \oeaep , e^{a\ep}(\uaep) ]\, dx + \qep 
\int_{\Omb} [A^b \oebep  , e^{b\ep}(\ubep) ]\, dx= \\ \\
\quad \quad =\int_{\Oma}\faep .\uaep  \, dx +\int_{\Omb}\fbep . \ubep\, dx +
 \int_{\Oma}[\gaep , e^{a\ep}(\uaep)]\, dx +
\int_{\Omb}[\gbep , e^{b\ep}(\ubep)]\, dx +\\ \\
\quad  \quad \quad\int_{\Sigma^a}\haep .\uaep  \, d\sigma  + 
\int_{\omb} \left( \hbepp .\ubep  _{ \;| x_3 =0}  
+ \hbepm . \ubep_{\;  | x_3 =-1} \right)\,dx'.
\end{array}
\end{equation}
Passing to the limit in the integral terms in $\Omb$ is easy. As for the terms in $\Oma \cap \{x_3>\rep\}$ and 
$\Sigma ^a\cap \{x_3>\rep\}$, we have from Lebesgue's theorem, with 
$\uaep=u^a+\rep v^a+(\rep)^2 w^a$ and 
$\chi^\ep$ the characteristic function of $\{x_3>\rep\}$, 
\[ \begin{array}{l}
\chi^\ep \eaep(\uaep) \rightarrow e^a(z^a) \mbox{ strongly in }
(L^2(\Oma))^{3\times 3}, \\  
\chi^\ep \uaep \rightarrow  u^a \mbox{ strongly in }(L^2(\Oma))^{3}, \\
\chi^\ep \uaep_{\; \;|\Sigma^a}\rightarrow  u^a _{\;|\Sigma^a}
\mbox{ strongly in }(L^2(\Sigma^a ))^{3},
\end{array}
\]
so that, by virtue of (\ref{convfa}), (\ref{convga}), (\ref{convha}), 
(\ref{defoe}),
\[\begin{array}{l}  
\int_{\Oma \cap \{x_3>\rep\}} [A^a  \oeaep-\gaep , 
e^{a\ep}(\uaep) ]\, dx 
-\int_{\Oma \cap \{x_3>\rep\}}\faep .\uaep  \, dx  
-\int_{\Sigma ^a\cap \{x_3>\rep\} } \haep.\uaep\, d\sigma  \\
\rightarrow \int_{\Oma}[A^a  \oea -g^a, 
e^a (z^a) ]\, dx -\int_{\Oma} f^a. u^a\, dx -\int_{\Sigma}h^a. u^a\, d\sigma.  
\end{array}\] 
For the terms in $\Oma \cap \{0<x_3<\rep\}$ and 
$\Sigma ^a\cap \{0<x_3<\rep\} $, it is clear, from (\ref{convfa}), 
(\ref{convha}) and from the uniform boundedness of $\uaep$, that 
\[  
\int_{\Oma \cap\{0<x_3<\rep\}  }  \faep .\uaep  \, dx 
-\int_{\Sigma ^a \cap\{0<x_3<\rep\}} \haep.\uaep\, d\sigma  
\rightarrow  0.
 \] 
Hence it remains to show that 
\[
\int_{\Oma \cap\{0<x_3<\rep\}  }[A^a  \oeaep-\gaep , 
e^{a\ep}(\uaep) ]\, dx \rightarrow 0.
\]
But we have, from Cauchy-Schwarz inequality, (\ref{convga}) and (\ref{defoe}),
\[
\int_{\Oma \cap\{0<x_3<\rep\}  }[A^a \oeaep -\gaep , 
e^{a\ep}(\uaep) ]\, dx \leq  C  \; \| e^{a\ep}(\uaep) \|_
  {(L^2(\Oma \cap\{0<x_3<\rep\} ))^{3 \times 3}}
\] 
and it is enough to prove that 
\begin{equation}\label{enough}
\| e^{a\ep}(\uaep) \|_
  {(L^2(\Oma\cap\{0<x_3<\rep\}))^{3 \times 3}}\rightarrow 0.
\end{equation}
Then, by passing to the limit in (\ref{test}), we will get
\begin{equation}\label{vg} \begin{array}{l}
\int_{\Oma} [A^a e^a(\oza), e^a(z^a)]\, dx+ 
 q \int_{\Omb} [A^b e^b(\ozb), e^b(z^b)]\, dx= \\ \\
\quad \quad 
\int_{\Oma}f^a . u^a\, dx +\int_{\Omb}f^b . u^b\, dx +
\int_{\Oma} [ g^a ,e^a(z^a) ]\, dx  +
\int_{\Omb} [ g^b, e^b(z^b)]\, dx  +\\ \\
\quad  \quad \quad\int_{\Sigma^a}h^a . u^a\, d\sigma + 
\int_{\omb} \left( h^b_+ . u^b_{ \;| x_3 =0}  
+ h ^b_- . u^b_{\;  | x_3 =-1} \right)\,dx',
\end{array}
\end{equation}
for any $z$ having the regularity given in (\ref{reg3}), and 
then, by a density argument given in Lemma 8 of the Appendix (Section 8.2), for any $z$ satisfying all 
the requirements of ${\cal Z}$, except the integral conditions. 
{\em A fortiori}, the same variational equality holds true for any $z$ in 
${\cal Z}$, that is $(\oza, \ozb)$ solves (\ref{vargrand1}).\\ \\

{\bf Proof of (\ref{enough}):}
We are going to prove that each term $e_{33}(\uaep)$, 
$\frac{1}{(\rep)^2}e_{\alpha \beta}(\uaep)$ and 
$\frac{1}{\rep}e_{\alpha 3}(\uaep)$ tend to zero (strongly) in 
$L^2(\{0<x_3<\rep\})$.\\

$\bullet$ {\bf Term $e_{33}(\uaep)$:}
We easily get from (\ref{uaep3})
\begin{equation}\label{e33}
e_{33}(\uaep) = \frac{\partial \uaep_3}{\partial x_3}= \frac{1}{\rep}
\left( u^a_3(x', \rep) - u^b_3 (\rep x')  \right) + v^a_3(x', \rep)
-\frac{(\ep)^2}{\rep}w^b_3 (\rep x',0).
\end{equation}
Clearly the last two terms in (\ref{e33}) tend to zero in 
$L^2(\{0<x_3<\rep\})$, since $v^a_3(x', \rep)$ and $w^b_3 (\rep x',0)$ 
are uniformly bounded:
\[
\int_{0 < x_3 <\rep} |v^a_3(x', \rep)|^2 \, dx \leq C \rep \rightarrow 0,
\]
\[
\int_{0 < x_3 <\rep} \left(\frac{(\ep)^2}{\rep}\right)^2
| w^b_3 (\rep x',0) |^2 \, dx \leq C \frac{\ep^4}{\rep} 
\ll C \frac{\ep^2}{\rep}\rightarrow 0,
\]
since, by assumption, $\ep^2 \ll \rep$. As for the first term in (\ref{e33}), 
it is uniformly bounded, because of the junction condition:
\[\begin{array}{l}
\frac{1}{\rep}
\left( u^a_3(x', \rep) - u^b_3 (\rep x')  \right) = \frac{1}{\rep}
\left( u^a_3(x', 0) + \int _0^{\rep} \frac{\partial u^a_3}{\partial x_3}
(x',t) \, dt - u^b_3 (0)-\int _0^{\rep} \nabla u^b_3 (t x').x'\, dt \right)\\ \\
=  
\frac{1}{\rep}\int _0^{\rep} \frac{\partial u^a_3}{\partial x_3}
(x',t) \, dt -\frac{1}{\rep}\int _0^{\rep} \nabla u^b_3 (t x').x'\, dt \leq C
\end{array}
\]
and hence, its norm in $L^2(\{0<x_3<\rep\})$ tends to zero.\\ 

$\bullet$ {\bf Term $\frac{1}{(\rep)^2}e_{\alpha \beta}(\uaep)$:}
From (\ref{uaep2}),
\[ \begin{array}{ll}
\frac{\partial \uaep _{\alpha}}{\partial x_{\beta}}&=
\frac{x_3}{\rep}
\left(\rep \frac{\partial v^a _{\alpha}}{\partial x_{\beta}}(x',\rep)
+(\rep)^2 \;\frac{\partial w^a _{\alpha}}{\partial x_{\beta}}(x',\rep)
\right) \\ \\
&+\left(1-\frac{x_3}{\rep}\right) \ep (\rep)^2
\left(\frac{\partial  \zeta ^b_{\alpha}}{\partial x_{\beta}}(\rep x') +
\ep\frac{\partial  v^b _{\alpha}}{\partial x_{\beta}}(\rep x',0)
\right)\end{array}
\]
and hence, since $e_{\alpha \beta}(v^a)=0$,
\[
\frac{1}{(\rep)^2}e_{\alpha \beta}(\uaep)=
\frac{x_3}{\rep}e_{\alpha \beta}(w^a)(x',\rep) + 
\left(1-\frac{x_3}{\rep}\right) \ep 
\left(e_{\alpha \beta}(\zeta^b)(\rep x')+ \ep e_{\alpha \beta}(v^b)(\rep x',0)
\right),
\]
which gives, from the regularity of $w^a$, $\zeta^b$ and $v^b$ (see (\ref{reg3}))
\[
\displaystyle \left|\frac{1}{(\rep)^2}e_{\alpha \beta}(\uaep)\right| 
\leq C + C \ep (1+\ep)\leq C,
\]
and hence, the norm of this term, in $L^2(\{0<x_3<\rep\})$, tends to zero.\\ 

$\bullet$ {\bf Term $\frac{1}{\rep}e_{\alpha 3}(\uaep)$:}
From (\ref{uaep2}),
\[  
\frac{\partial \uaep_{\alpha}}{\partial x_3}=\frac{1}{\rep}
\left(
u^a_{\alpha}(\rep) + \rep v^a_{\alpha}(x',\rep) +(\rep)^2 w^a_{\alpha}(x',\rep)
\right) - \ep  \left(\zeta^b_{\alpha}(\rep x')+ \ep v^b_{\alpha}(\rep x',0)
\right)
\]
and, from (\ref{uaep3}),
\[ 
\frac{\partial \uaep_{3}}{\partial x_{\alpha}}=\frac{x_3}{\rep}
\left(\frac{\partial u^a_{3}}{\partial x_{\alpha}}(x',\rep)+\rep
\frac{\partial v^a_{3}}{\partial x_{\alpha}}(x',\rep)
\right) + \left(1-\frac{x_3}{\rep}\right)\rep
\left(\frac{\partial u^b_{3}}{\partial x_{\alpha}}(\rep x')+\ep^2
\frac{\partial w^b_{3}}{\partial x_{\alpha}}(\rep x',0)
\right),
\]
so that 
\begin{equation}\label{ea3}
\frac{2}{\rep}e_{\alpha 3}(\uaep)=\frac{1}{\rep}\left(
\frac{\partial \uaep_{\alpha}}{\partial x_3} +
\frac{\partial \uaep_{3}}{\partial x_{\alpha}}
\right)=T_1 +T_2 +T_3 +T_4,
\end{equation}
with
\[\begin{array}{l}
T_1= \frac{1}{(\rep)^2}\left(u^a_{\alpha}(\rep)+x_3
\frac{\partial u^a_{3}}{\partial x_{\alpha}}(x',\rep)\right) ,\\ \\
T_2=\frac{1}{\rep}\left(v^a_{\alpha}(x',\rep)-\ep 
\left(\zeta^b_{\alpha}(\rep x') +\ep v^b_{\alpha}(\rep x',0)  \right) \right),\\ \\
T_3=w^a_{\alpha}(x',\rep),\\ \\
T_4=\frac{x_3}{\rep}\frac{\partial u^a_{3}}{\partial x_{\alpha}}(x',\rep)+
\left(1-\frac{x_3}{\rep} \right)
\left(
\frac{\partial u^b_{3}}{\partial x_{\alpha}}(\rep x')+\ep^2
\frac{\partial w^b_{3}}{\partial x_{\alpha}}(\rep x',0)  
\right).
\end{array}
\]
We are going to show that each term tends to zero, in $L^2(\{0<x_3<\rep\})$.\\ 

$\circ$ {\bf Term $T_1$:} As $u^a_{\alpha}(0)=0$,
\[  
u^a_{\alpha}(\rep)=\int_{0}^{\rep}\frac{d u^a_{\alpha}}{dx_3}(t)\, dt
\]
and, as $e_{\alpha 3}(u^a)=0$,
\[  
x_3\frac{\partial u^a_{3}}{\partial x_{\alpha}}(x',\rep)=-x_3
\frac{d u^a_{\alpha}}{dx_3}(\rep)=-\frac{x_3}{\rep}\int_{0}^{\rep}
\frac{d u^a_{\alpha}}{dx_3}(t)\,dt +\frac{x_3}{\rep}\int_{0}^{\rep}
\left(
\frac{d u^a_{\alpha}}{dx_3}(t) -\frac{d u^a_{\alpha}}{dx_3}(\rep) 
\right)\, dt,
\]
so that, since $\frac{d u^a_{\alpha}}{dx_3}(0)=0$,
\[\begin{array}{ll}
u^a_{\alpha}(\rep)+x_3\frac{\partial u^a_{3}}{\partial x_{\alpha}}(x',\rep)
&=\left(1- \frac{x_3}{\rep}\right)\displaystyle \int_{0}^{\rep}
\frac{d u^a_{\alpha}}{dx_3}(t)\,dt +\frac{x_3}{\rep}\displaystyle\int_{0}^{\rep}
\left(
\frac{d u^a_{\alpha}}{dx_3}(t) -\frac{d u^a_{\alpha}}{dx_3}(\rep) 
\right)\, dt \\ \\
&= \displaystyle\int_{0}^{\rep}\int_{0}^{t}
\frac{d^2 u^a_{\alpha}}{dx_3^2}(\tau)\, d\tau \, dt
\end{array}
\]
and, from the regularity of $u^a_{\alpha}$,
\[
\left|T_1\right| =\frac{1}{(\rep)^2}
\left|
u^a_{\alpha}(\rep)+x_3\frac{\partial u^a_{3}}{\partial x_{\alpha}}(x',\rep)  
\right| \leq \frac{1}{(\rep)^2}\displaystyle \int_{0}^{\rep}\int_{0}^{\rep}
\left|
\frac{d^2 u^a_{\alpha}}{dx_3^2}(\tau)
\right|\, d\tau \, dt \leq C, 
\]
so that $T_1$ tends to zero, in $L^2(\{0<x_3<\rep\})$.\\

$\circ$ {\bf Term $T_2$:} We have
\[
T_2 =\frac{1}{\rep}v^a_{\alpha}(x',\rep) -\frac{\ep}{\rep}
\left(\zeta^b_{\alpha}(\rep x') +\ep v^b_{\alpha}(\rep x',0)  \right).
\]
But, as $c(0)=0$,
\[
\left|\frac{1}{\rep}v^a_{\alpha}(x',\rep) \right| 
\leq  \frac{C}{\rep} \left|c(\rep)\right|\leq C
\]
and this term tends to zero, in $L^2(\{0<x_3<\rep\})$. 
Moreover, as $\ep^2 \ll \rep$ and as $\zeta^b_{\alpha}$ and $v^b_{\alpha}$ 
are uniformly bounded, due to the regularity conditions (\ref{reg3}),
\[
\left| \frac{\ep}{\rep}
\left(\zeta^b_{\alpha}(\rep x') +\ep v^b_{\alpha}(\rep x',0)  \right) \right|
\leq C \frac{\ep}{\rep},
\]
so that the $L^2(\{0<x_3<\rep\})$-norm of this term is bounded by $ C\ep /
\sqrt{\rep}$, which tends to zero by assumption.\\
 
$\circ$ {\bf Terms $T_3$ and $T_4$:} Clearly these terms are bounded, 
due to the regularity conditions (\ref{reg3}).

%\sub
\section{Proof of stronger convergences and proof of Corollary 1}

Actually, the stronger convergences in Theorem 1 are deduced from Corollary 1. 
The proof preceeds as follows.

Taking $\ouep=(\ouaep,\oubep)$ as test function in the variational equation 
of Problem (\ref{vargrand}), we get
\begin{equation}\label{energy} \begin{array}{l}
{\cal E}^\ep = \int_{\Oma} [A^a \oeaep , \oeaep  ]\, dx + \qep 
\int_{\Omb} [A^b \oebep  ,  \oebep ]\, dx= \\ \\
\quad \quad =\int_{\Oma}\faep .\ouaep  \, dx +\int_{\Omb}\fbep . \oubep\, dx +
 \int_{\Oma}[\gaep , \oeaep ]\, dx +
\int_{\Omb}[\gbep , \oebep ]\, dx +\\ \\
\quad  \quad \quad\int_{\Sigma^a}\haep .\ouaep  \, d\sigma  + 
\int_{\omb} \left( \hbepp .\oubep  _{ \;| x_3 =0}  
+ \hbepm . \oubep_{\;  | x_3 =-1} \right)\,dx'.
\end{array}
\end{equation}
We are going to pass to the limit in the last member of the above equality.\\

$\bullet$ If $q \in (0,+\infty)$, we have, from the convergences already 
proved in Theorem 1 and from classical compactness arguments,
\[\begin{array}{l}
(\oeaep,\oebep) \rightharpoonup (\oea,\oeb) \quad\mbox{ weakly in }
(L^2(\Oma))^{3 \times 3} \times (L^2(\Omb))^{3 \times 3}, \\ \\
(\ouaep,\oubep) \rightarrow (\oua,\oub) \quad\mbox{ strongly in }
(L^2(\Oma))^{3} \times (L^2(\Omb))^{3}, \\ \\
\ouaep_{\; \;| \Sigma ^a}\rightarrow \oua _{\; | \Sigma ^a}\quad\mbox{ strongly in }
(L^2(\Sigma ^a))^{3},\\ \\
\oubep_{\;\;| x_3 =0} \quad(\mbox{ resp }\oubep_{\;\;| x_3 =-1})\rightarrow
\oub _{\; |x_3 =0} \quad(\mbox{ resp }\oub _{\; |x_3 =-1})\quad\mbox{ strongly in }
(L^2(\omb))^{3}.
\end{array}
\]
If $(\gaep,\gbep)$ tends to $(g^a,g^b)$ strongly in 
$(L^2(\Oma))^{3 \times 3} \times (L^2(\Omb))^{3 \times 3}$, it follows that
\[\begin{array}{l}
{\cal E}^\ep=\int_{\Oma}\faep .\ouaep  \, dx +\int_{\Omb}\fbep . \oubep\, dx +
 \int_{\Oma}[\gaep , \oeaep ]\, dx +
\int_{\Omb}[\gbep , \oebep ]\, dx + 
 \int_{\Sigma^a}\haep .\ouaep  \, d\sigma  + \\ \\
\qquad\qquad\qquad +\int_{\omb} \left( \hbepp .\oubep  _{ \;| x_3 =0}  
 + \hbepm . \oubep_{\;  | x_3 =-1} \right)\,dx' \longrightarrow \\ \\
\int_{\Oma} f^a .\oua  \, dx +\int_{\Omb} f^b . \oub \, dx +
 \int_{\Oma}[ g^a , \oea  ]\, dx +
\int_{\Omb}[g^b , \oeb ]\, dx + 
 \int_{\Sigma^a}h^a .\oua   \, d\sigma  + \\ \\
\qquad\qquad\qquad + \int_{\omb} \left(  h^b_+ .\oub   _{ \;+| x_3 =0}  
+  h^b_- . \oub _{\;  _| x_3 =-1} \right)\,dx'=\\ \\
 = \int_{\Oma} [A^a \oea  , \oea  ]\, dx + q 
\int_{\Omb} [A^b \oeb  ,  \oeb ]\, dx={\cal E},
\end{array}
\]
which proves the first part of Corollary 1. Moreover, we get, from the 
convergence of ${\cal E}^\ep$ to ${\cal E}$ and from a classical 
lower semicontinuity argument:
\[\begin{array}{ll}
0&=\liminf \left(\int_{\Oma} [A^a \oeaep  ,  \oeaep ]\, dx   -
\int_{\Oma} [A^a \oea  , \oea  ]\, dx +
\qep \int_{\Omb} [A^b \oebep  ,  \oebep ]\, dx -
q \int_{\Omb} [A^b \oeb  ,  \oeb ]\, dx
\right)
\\ \\
&\geq \liminf\left(\int_{\Oma} [A^a \oeaep  ,  \oeaep ]\, dx  -
\int_{\Oma} [A^a \oea  , \oea  ]\, dx \right) + \\ \\
 &\qquad\qquad\liminf\left(\qep \int_{\Omb} [A^b \oebep  ,  \oebep ]\, dx -
q \int_{\Omb} [A^b \oeb  ,  \oeb ]\, dx
\right)\\ \\
&=\liminf\left(\int_{\Oma} [A^a \oeaep  ,  \oeaep ]\, dx  -
\int_{\Oma} [A^a \oea  , \oea  ]\, dx \right) + \\ \\
&\qquad\qquad\liminf q\left(  \int_{\Omb} [A^b \oebep  ,  \oebep ]\, dx -
 \int_{\Omb} [A^b \oeb  ,  \oeb ]\, dx
\right) \geq 0,
\end{array} 
\]
which gives, up to extraction of a new subsequence,
\[\begin{array}{l}
  \int_{\Oma} [A^a \oeaep  ,  \oeaep ]\, dx  \longrightarrow
\int_{\Oma} [A^a \oea  , \oea  ]\, dx  , \\ \\ 
   \int_{\Omb} [A^b \oebep  ,  \oebep ]\, dx \longrightarrow 
 \int_{\Omb} [A^b \oeb  ,  \oeb ]\, dx.
\end{array} 
\]
It follows that (e.g.)
\[\begin{array}{l}
 C \; \| \oeaep -\oea \|^2 _{(L^2(\Oma))^{3 \times 3}}
\leq \int_{\Oma} [A^a (\oeaep -\oea)  , (\oeaep -\oea) ]\, dx= \\ \\
=\int_{\Oma} [A^a \oeaep  , \oeaep  ]\, dx + 
\int_{\Oma} [A^a \oea  , \oea  ]\, dx 
-\int_{\Oma} [A^a \oeaep  , \oea  ]\, dx -
\int_{\Oma} [A^a \oea , \oeaep  ]\, dx\longrightarrow 0,
\end{array} 
\]
and hence $\oeaep$ tends to $\oea$ strongly in $(L^2(\Oma))^{3 \times 3}$. 
It follows that $e(\ouaep) \rightarrow e(\oua)$ strongly in 
$(L^2(\Oma))^{3 \times 3}$ and then, from Korn's inequality, $\ouaep 
\rightarrow \oua$ strongly in $H^1(\Oma)^{3}$. 
By the same way, $\oebep$ tends to $\oeb$ strongly in 
$(L^2(\Omb))^{3 \times 3}$ and $\oubep 
\rightarrow \oub$ strongly in $H^1(\Omb)^{3}$. The conclusion is that we get the stronger 
convergences mentionned in Theorem 1, if $q \in (0,+\infty)$.\\

$\bullet$ If $q =+\infty$, we have seen that
\[
\oubep \rightarrow \oub=0 \mbox{ strongly in } (H^1(\Omb))^3, 
\] 
\[
\oebep \rightarrow \oeb = 0 \mbox{ strongly in } (L^2(\Omb))^{3 \times 3}, 
\]
and, with appropriate changes in the above proof, we have, if $\gaep$ tend to 
$g^a$ strongly in $(L^2(\Oma))^{3 \times 3}$,
\[ 
{\cal E}^\ep \longrightarrow  
\int_{\Oma} f^a .\oua  \, dx  +
 \int_{\Oma}[ g^a , \oea  ]\, dx +
 \int_{\Sigma^a}h^a .\oua   \, d\sigma   
 = \int_{\Oma} [A^a \oea  , \oea  ]\, dx   ={\cal E}_{\infty},
\]
\[\begin{array}{ll}
0&=\liminf \left(\int_{\Oma} [A^a \oeaep  ,  \oeaep ]\, dx   -
\int_{\Oma} [A^a \oea  , \oea  ]\, dx +
\qep \int_{\Omb} [A^b \oebep  ,  \oebep ]\, dx  
\right)
\\ \\
&\geq \liminf\left(\int_{\Oma} [A^a \oeaep  ,  \oeaep ]\, dx  -
\int_{\Oma} [A^a \oea  , \oea  ]\, dx \right) + 
\liminf\left(\qep \int_{\Omb} [A^b \oebep  ,  \oebep ]\, dx  
\right)  \geq 0,
\end{array} 
\]
\[\begin{array}{l}
  \int_{\Oma} [A^a \oeaep  ,  \oeaep ]\, dx  \longrightarrow
\int_{\Oma} [A^a \oea  , \oea  ]\, dx  , \\ \\ 
 \qep  \int_{\Omb} [A^b \oebep  ,  \oebep ]\, dx \longrightarrow 0, \\ \\
 \oeaep \rightarrow \oea \mbox{ strongly in }(L^2(\Oma))^{3 \times 3}, \\ \\
\sqrt{\qep}\oebep \rightarrow 0\mbox{ strongly in }(L^2(\Oma))^{3 \times 3}, 
\\ \\
\ouaep \rightarrow \oua \mbox{ strongly in }H^1(\Oma)^3.
\end{array} 
\]\\

$\bullet$ If $q =0$, we have, with $\utep =\qep \ouep$, $\etep =\qep \oeep$,
\begin{equation}\label{energy} \begin{array}{l}
{\cal E}^\ep = \frac{1}{\qep}\int_{\Oma} [A^a \etaep , \etaep  ]\, dx +  
\int_{\Omb} [A^b \etbep  ,  \etbep ]\, dx= \\ \\
\quad \quad =\int_{\Oma}\faep .\utaep  \, dx +\int_{\Omb}\fbep . \utbep\, dx +
 \int_{\Oma}[\gaep , \etaep ]\, dx +
\int_{\Omb}[\gbep , \etbep ]\, dx +\\ \\
\quad  \quad \quad\int_{\Sigma^a}\haep .\utaep  \, d\sigma  + 
\int_{\omb} \left( \hbepp .\utbep  _{ \;| x_3 =0}  
+ \hbepm . \utbep_{\;  | x_3 =-1} \right)\,dx',
\end{array}
\end{equation}
\[
\utaep \rightarrow \oua=0 \mbox{ strongly in } (H^1(\Oma))^3, 
\] 
\[
\etaep \rightarrow \oea = 0 \mbox{ strongly in } (L^2(\Oma))^{3 \times 3}, 
\]
and  we have, if $\gbep$ tend to 
$g^b$ strongly in $(L^2(\Omb))^{3 \times 3}$,
\[ 
{\cal E}^\ep \longrightarrow  
\int_{\Omb} f^b .\oub  \, dx  +
 \int_{\Omb}[ g^b , \oeb  ]\, dx +
\int_{\omb} \left( h^b_+ .\oub  _{ \;| x_3 =0}  
+  h^b_-  . \oub_{\;  | x_3 =-1} \right)\,dx'  
 = \int_{\Omb} [A^b \oeb  , \oeb  ]\, dx   ={\cal E}_0,
\]
\[\begin{array}{ll}
0&=\liminf \left( \frac{1}{\qep } \int_{\Oma} [A^a \etaep  , \etaep  ]\, dx + 
\int_{\Omb} [A^b \etbep   ,  \etbep ]\, dx -  
\int_{\Omb} [A^b \oeb  , \oeb  ]\, dx  
\right)
\\ \\
&\geq \liminf\left(\frac{1}{\qep } \int_{\Oma} [A^a \etaep  , \etaep  ]\, dx 
\right) + 
\liminf\left( \int_{\Omb} [A^b \etbep   ,  \etbep ]\, dx -  
\int_{\Omb} [A^b \oeb  , \oeb  ]\, dx  
\right)  \geq 0,
\end{array} 
\]
\[\begin{array}{l}
  \int_{\Omb} [A^b \etbep  ,  \etbep  ]\, dx  \longrightarrow
\int_{\Omb} [A^a \oeb  , \oeb  ]\, dx  , \\ \\ 
\frac{1}{\qep}  \int_{\Oma} [A^a \etaep, \etaep ]\, dx \longrightarrow 0, \\ \\
\qep \oebep = \etbep\rightarrow \oeb\mbox{ strongly in }(L^2(\Omb))^{3 \times 3}, \\ \\
\sqrt{\qep}\;\oeaep=\frac{1}{\sqrt{\qep}} \etaep\rightarrow 0\mbox{ strongly in }
(L^2(\Omb))^{3 \times 3},\\ \\
\qep \oubep \rightarrow \oub \mbox{ strongly in }H^1(\Omb)^3.
\end{array} 
\]\\

\section{Appendix}
 
\subsection{The definitions of $(v^a, w^a)$ and  $(v^b, w^b)$ 
as suitable limits}

For the convenience of the reader, we give in this Appendix a sketch of proof of the following result, 
mentionned in Section 4.2. (For thin cylinders, a complete proof can be find in 
\cite{MuSi2}. The case of plates is analogous and simpler, but it is not published, 
as far as we know.)

\begin{lemma} (i) Let $\{u^\ep\}_\ep$ be a sequence in $(H^1(\Oma))^3$, such that
$u^\ep =0$ on $T^a=\oma \times\{1\}$ and
\begin{equation}\label{91}
\{\eaep(u^\ep)\}_\ep \mbox{ is bounded in }(L^2(\Oma))^{3 \times 3}.
\end{equation}
Let ${\cal W}^a$ be the space defined in Section 2.2 and let
\[
\begin{array}{l} 
{\cal V}^a_-= 
\{
 v^a \in (H^1(\Oma))^2 \times L^2(0,1;H^1 (\oma)), \; 
\exists c \in H^1(0,1), c(1)=0,
v^a_1=-c\,x_2, v^a_2=c\,x_1,  \\
\qquad\qquad\qquad\qquad \mbox{ for a.e. } 
x_3 \in (0,1), \; 
 \int_{\oma}v^a_3 (x',x_3)\, dx'=0 
\}.
\end{array} 
\]
(Note that ${\cal V}^a_-$ satisfies the same requirements as ${\cal V}^a$, 
in Section 2.2, except $c(0)=0$.)
Then there exists a pair $(v^a,w^a) \in {\cal V}^a_- \times {\cal W}^a$, 
such that, for all $\alpha, \beta=1,2$:
\begin{equation}\label{92}
\frac{1}{\rep}e_{\alpha 3}(u^\ep) \rightharpoonup e_{\alpha 3}(v^a) 
\mbox{ weakly in }L^2(\Oma), 
\end{equation}
\begin{equation}\label{93}
\frac{1}{(\rep)^2}e_{\alpha \beta}(u^\ep) \rightharpoonup e_{\alpha \beta}(w^a) 
\mbox{ weakly in }L^2(\Oma).
\end{equation}
%In addition, there exist $\tilde{u}^\ep$ and $\hat{u}^\ep$, such that, 
%for all $\alpha, \beta=1,2$:
%\begin{equation}\label{94}
%e_{\alpha 3}(\tilde{u}^\ep)=0, \; e_{\alpha \beta}(\hat{u}^\ep)=0 
%\end{equation}
%and 
%\begin{equation}\label{95}
%\left\{ 
%\begin{array}{l}
%\frac{u^\ep_\alpha}{\rep} - \tilde{u}^\ep_\alpha
%- \displaystyle\int_0^1 \left(
% \frac{u^\ep_\alpha}{\rep} - \tilde{u}^\ep_\alpha \right)\,dx_3
%\rightharpoonup l_\alpha
%\mbox{ weakly in }L^2(\oma;H^1(0,1)), \\ \\
%\frac{\partial l_\alpha}{\partial x_3}=2 e_{\alpha 3}(v^a),
%\end{array}
%\right.
%\end{equation}
Morever, denoting by $(c, v^a_3)$ the couple defining $v^a$ and setting 
\begin{equation}\label{cep}
c^\ep(x_3) =  \frac
{\int_{\oma}\left(x_1  u^\ep_2(x',x_3) 
-x_2  u^\ep_1(x',x_3) \right)\, dx'}
{\rep\int_{\oma}\left(x_1^2 + x_2^2 \right)\, dx'},
\end{equation}
\begin{equation}\label{vaep}
v^\ep_3 = \frac{ u^\ep_3}{\rep} -
\frac{1}{|\oma|}\int_{\oma}\frac{ u^\ep_3}{\rep} \, dx'+
\frac{1}{|\oma|}
\sum _{\alpha}x_{\alpha } \frac{d}{dx_3}  \int_{\oma}u^\ep_\alpha \, dx',
\end{equation}
we have 
\begin{equation}\label{c}
c^\ep \rightarrow c \mbox{ strongly in } L^2(0,1),
\end{equation}
\begin{equation}\label{va3}
v^\ep_3 \rightharpoonup v^a_3 \mbox{ weakly in } H^{-1}(0,1;H^1(\oma)).
\end{equation}
Finally, setting
\begin{equation}\label{dep}
d^\ep_\alpha(x_3) = \frac{1}{|\oma|}\int_{\oma} 
\frac{u^\ep_\alpha (x',x_3)}{\rep}\, dx'
\end{equation}
and $x^R_1=-x_2$, $x^R_2=x_1$, we have
\begin{equation}\label{96}
 \frac{u^\ep _\alpha}{(\rep)^2} - 
\frac{1}{\rep} \left( c^\ep x^R_\alpha + d^\ep_\alpha\right)
\rightharpoonup w^a_\alpha
\mbox{ weakly in }L^2(0,1;H^1(\oma)). 
\end{equation} 
(ii) If $\{u^\ep\}_\ep$ is a sequence in $(H^1(\Omb))^3$, such that
$u^\ep =0$ on $\Sigma^b=\partial \omb \times (-1,0)$ and
\begin{equation}\label{97}
\{\ebep(u^\ep)\}_\ep \mbox{ is bounded in }(L^2(\Omb))^{3 \times 3},
\end{equation}
then there exists a pair $(v^b,w^b) \in {\cal V}^b \times {\cal W}^b$, 
such that, for all $\alpha=1,2$:
\begin{equation}\label{98}
\frac{1}{\ep}e_{\alpha 3}(u^\ep) \rightharpoonup e_{\alpha 3}(v^b) 
\mbox{ weakly in }L^2(\Omb), 
\end{equation}
\begin{equation}\label{99}
\frac{1}{\ep^2}e_{33}(u^\ep) \rightharpoonup e_{33}(w^b) 
\mbox{ weakly in }L^2(\Omb).
\end{equation}
In addition, we have
\begin{equation}\label{911}
\frac{u^\ep _\alpha }{\ep} - \tilde{u}^\ep  _\alpha 
- \int_{-1}^{0}\left(\frac{u^\ep_\alpha}{\ep} - \tilde{u}^\ep_\alpha \right)
  \,dx_3
\rightharpoonup v^b_\alpha
\mbox{ weakly in }L^2(\omb;H^1(-1,0)), \; \mbox{ for } \alpha =1,2,
\end{equation}
with $\tilde{u}^\ep$ defined by 
\[
\tilde{u}^\ep_{\alpha}=
-\int_0^{x_3} \frac{1}{\ep} \frac{\partial u^\ep_3}{\partial x_\alpha}(x',s)ds.
\]
%\begin{equation}\label{911}
%\left\{ 
%\begin{array}{l}
% \frac{u^\ep _\alpha }{\ep} - \tilde{u}^\ep  _\alpha 
%- \int_{-1}^{0}\left(\frac{u^\ep_\alpha}{\ep} - \tilde{u}^\ep_\alpha \right)
%  \,dx_3
%\rightharpoonup t_\alpha
%\mbox{ weakly in }L^2(\omb;H^1(-1,0)), \\ \\
%\frac{\partial t_\alpha}{\partial x_3}=2 e_{\alpha 3}(v^b),
%\end{array}
%\right.
%\end{equation}
Moreover,
\begin{equation}\label{910}
\frac{u^\ep_3}{\ep^2} - \int_{-1}^{0} \frac{u^\ep_3}{\ep^2}\,dx_3
\rightharpoonup w^b_3
\mbox{ weakly in }L^2(\omb;H^1(-1,0)). 
\end{equation}
\end{lemma} 

{\em Proof of (i):} We use the following decomposition and estimate, 
whose proof may be found for instance in \cite{Le1} or in \cite{Le2}: 
there exists a positive constant $C$ such that, 
for every $u$ in $L^2(0,1;H^1(\oma))^2$, there exist $\overline{u}$ and 
$\hat{u}$ satisfying:
\begin{equation}\label{912}
\left\{
\begin{array}{l}
u=\overline{u}+\hat{u}, \; \\
\int_{\oma}\overline{u}_\alpha(x',x_3)\, dx'\equiv 0, \; 
\int_{\oma}\left(x_1\overline{u}_2(x',x_3)-
x_2\overline{u}_1(x',x_3)\right)\, dx'\equiv 0,  \\ 
e_{\alpha \beta}(\hat{u})=0, \qquad\forall \alpha, \beta =1,2,
\end{array}
\right.
\end{equation}
\begin{equation}\label{913}
\|\overline{u}\|_{(L^2(0,1;H^1(\oma)))^2} \leq C \sum_{\alpha, \beta}
\|e_{\alpha \beta}(u)\|_{L^2(\Oma)}.
\end{equation} 
The function $\hat{u}$ is a rigid displacement:
\begin{equation}\label{914}
\hat{u}_\alpha (x',x_3)= c(x_3) x^R_\alpha + d_\alpha (x_3),
\end{equation}
with $x^R_1=-x_2$, $x^R_2=x_1$ ($R$ for "rotation"). Applying (\ref{912}) and 
(\ref{914}) to $u=\frac{1}{\rep}(u^\ep_1,u^\ep_2)$, we get: 
\begin{equation}\label{915}
\frac{1}{\rep}u^\ep_\alpha = \overline{u}^\ep_\alpha + \hat{u}^\ep_\alpha, \mbox{ with }
\hat{u}^\ep_\alpha =c^\ep(x_3) x^R_\alpha + d^\ep_\alpha (x_3).
\end{equation}
One can check easily that the functions $c^\ep$ and $d^\ep_\alpha$ are given 
in terms of 
$u^\ep$ by the formulae (\ref{cep}) and (\ref{dep}).
%\begin{equation}\label{916}
%c^\ep(x_3) =  \frac
%{\int_{\oma}\left(x_1  u^\ep_2(x',x_3) 
%-x_2  u^\ep_1(x',x_3) \right)\, dx'}
%{\rep\int_{\oma}\left(x_1^2 + x_2^2 \right)\, dx'}, \;
%d^\ep_\alpha(x_3) = \frac{1}{|\oma|}\int_{\oma} 
%\frac{u^\ep_\alpha (x',x_3)}{\rep}\, dx'.
%\end{equation}
From (\ref{913}), we obtain:
\[
%\begin{equation}\label{917}
\|\overline{u}^\ep_\alpha\|_{L^2(0,1;H^1(\oma))} \leq C \sum_{\alpha, \beta}
\|e_{\alpha \beta}(\frac{1}{\rep}u^\ep)\|_{L^2(\Oma)}. 
%\end{equation}
\]
Setting $w^\ep_\alpha = \overline{u}^\ep_\alpha /\rep$ and using (\ref{91}), 
it follows that
\[
%\begin{equation}\label{919}
\|w^\ep_\alpha\|_{L^2(0,1;H^1(\oma))} \leq C. 
%\end{equation}
\]
So, taking a subsequence of $\ep$, still denoted by the same letter, 
we may assume the existence of $w^a_\alpha$ such that
\[
%\begin{equation}\label{920}
w^\ep_\alpha \rightharpoonup w^a_\alpha \mbox{ weakly in }L^2(0,1;H^1(\oma)), 
\forall \alpha=1,2,
%\end{equation}
\]
that is  (\ref{96}). Moreover it is clear that $(w^\ep_1,w^\ep_2,0)$ and 
$w^a = (w^a_1,w^a_2,0)$ belong to ${\cal W}^a$. Since (\ref{915}) implies that 
\[
\frac{1}{(\rep)^2}e_{\alpha \beta}(u^\ep)=e_{\alpha \beta}(w^\ep),
\]
we see that (\ref{93}) is proved.  

It remains to prove the convergences involving $v^a$. In Section 5.3, it is proved that there exists $c$ in $H^1(0,1)$, $c(1)=0$, such that, 
for a subsequence of $\ep$, (\ref{c}) holds true. As for the other convergences 
involving $v^a$, we use again the decomposition (\ref{915}), 
from which we deduce the following equality: 
\begin{equation}\label{921}
\frac{2}{\rep}e_{\alpha 3}(u^\ep) =  
  \frac{\partial \overline{u}^\ep_\alpha}{\partial x_3} +
\frac{d c^\ep}{dx_3}x^R_\alpha +\frac{d d^\ep_\alpha}{dx_3}
  +\frac{1}{\rep} \frac{\partial  u^\ep_3}{\partial x_\alpha},
 \forall \alpha=1,2.  
\end{equation}
Now, setting
\[  
v^\ep_3 = \frac{ u^\ep_3}{\rep} -
\frac{1}{|\oma|}\int_{\oma}\frac{ u^\ep_3}{\rep} \, dx'+
x_{\beta} \frac{d}{dx_3} d^\ep_{\beta}(x_3)
\]
(the summation convention is used, concerning the index $\beta$),
equality (\ref{921}) can be written as:
\begin{equation}\label{927}
\frac{2}{\rep}e_{\alpha 3}(u^\ep)= \frac{d c^\ep}{dx_3}x^R_\alpha+
\frac{\partial v^\ep_3}{\partial x_\alpha}
 +   \frac{\partial \overline{u}^\ep_\alpha}{\partial x_3}.
\end{equation}
The following estimate is proved in \cite{MoMuSi1}:
\[
%\begin{equation}\label{928}
\| v^\ep_3\|_{H^{-1}(0,1;H^1(\oma))} \leq  C
\left(\sum_{\alpha \beta} \|e_{\alpha \beta}(\frac{u^\ep}{\rep})\|_{L^2(\Oma)}+
 \sum_{\alpha }\|e_{\alpha 3}(\frac{u^\ep}{\rep})\|_{L^2(\Oma)}
\right),
%\end{equation}
\]
so that, from (\ref{91}), the sequence 
$\{v^\ep_3\}_\ep$ is bounded in $H^{-1}(0,1;H^1(\oma))$. Hence there exists 
$v^a_3$ in $H^{-1}(0,1;H^1(\oma))$, having zero mean-value on $\oma$, 
such that (\ref{va3}) holds true (for a subsequence).
%\begin{equation}\label{929}
%v^\ep_3 \rightharpoonup v^a_3 \mbox{ weakly in }H^{-1}(0,1;H^1(\oma)),
%\end{equation}
It follows also from (\ref{91}) that  
\begin{equation}\label{929'}
\frac{1}{\rep}e_{\alpha 3}(u^\ep) \rightharpoonup \tau_{\alpha 3} 
\mbox{ weakly in } L^2(\Oma)
\end{equation}
(again for some subsequence and 
some $\tau_{\alpha 3}$ in $L^2(\Oma)$).
Moreover, since $w^\ep_\alpha$ is bounded in $L^2(0,1;H^1(\oma))$,  
\begin{equation}\label{929"} 
 \frac{\partial \overline{u}^\ep_\alpha}{\partial x_3}=\rep
\frac{\partial  w^\ep_\alpha}{\partial x_3}
\mbox { tends to }0
\mbox { in the sense of distributions }.
\end{equation} 
By passing to the limit in (\ref{927}), using (\ref{c}), (\ref{va3}), 
(\ref{929'}) and (\ref{929"}), 
we get:
\begin{equation}\label{930}
2 \tau_{\alpha 3} = \frac{dc}{dx_3}x^R_\alpha + 
\frac{\partial}{\partial x_\alpha} v^a_3,
\end{equation}
which implies that
\begin{equation}\label{931}
\frac{\partial}{\partial x_\alpha} v^a_3 \in L^2(\Oma).
\end{equation}
From (\ref{931}) and from the fact that $v^a_3$ is in $H^{-1}(0,1;H^1(\oma))$, 
it is then an exercise to get that $v^a_3$ is in $L^2(0,1;H^1(\oma))$, so that 
$v^a:=(c(x_3)x^R_\alpha, v^a_3) \in {\cal V}^a_-$ and satisfies (\ref{92}).  

%It remains to prove (\ref{95}). Defining $\tilde{u}^\ep$ by 
%\begin{equation}\label{utilde}
%\tilde{u}^\ep=\left(
%-\int_0^{x_3} \frac{1}{\rep} \frac{\partial u^\ep_3}{\partial x_1}(x',s)ds,
%-\int_0^{x_3} \frac{1}{\rep} \frac{\partial u^\ep_3}{\partial x_2}(x',s)ds,
%\frac{u^\ep_3}{\rep}
%\right),
%\end{equation}
%it obviously satisfies
%\begin{equation}\label{932}
%e_{\alpha 3}(\tilde{u}^\ep)=0,\;\forall \alpha=1,2.
%\end{equation}
%The sequence of functions 
%\begin{equation}\label{l}
%l^\ep_{\alpha}=
%\frac{u^\ep_{\alpha}}{\rep}- \tilde{u}^\ep_{\alpha}-\int_0^1\left(
%\frac{u^\ep_{\alpha}}{\rep} - \tilde{u}^\ep_{\alpha}\right)dx_3
%\end{equation}
%is bounded in 
%$L^2(\oma;H^1(0,1))$, since it has mean value zero with respect to $x_3$ 
%and since
%\[
%\frac{\partial l^\ep_{\alpha}}{\partial x_3}=\frac{2}{\rep}e_{\alpha 3}(u^\ep)
%\] 
%is bounded in $L^2(\Oma)$. Taking into account (\ref{92}), we get (\ref{95}). 

{\em Proof of (ii):} Now we prove the analogous of the previous properties 
in the framework of 3d-2d reduction of dimension. This is much easier. Indeed, 
in order to prove 
(\ref{98}) and (\ref{911}), we consider the sequence $\{v^\ep_\alpha\}_\ep$ 
defined by:
\[
%\begin{equation}\label{934}
v^\ep_\alpha =\frac{u^\ep_\alpha}{\ep}-\tilde{u}^\ep_\alpha - 
\int_{-1}^0 \left(\frac{u^\ep_\alpha}{\ep}-\tilde{u}^\ep_\alpha\right)
%\end{equation}
\]
and 
\[
%\begin{equation}\label{utilde'}
\tilde{u}^\ep=\left(
-\int_0^{x_3} \frac{1}{\ep} \frac{\partial u^\ep_3}{\partial x_1}(x',s)ds,
-\int_0^{x_3} \frac{1}{\ep} \frac{\partial u^\ep_3}{\partial x_2}(x',s)ds,
\frac{u^\ep_3}{\ep}
\right),
%\end{equation}
\]
Then we have as above that
\begin{equation}\label{935}
\frac{\partial v^\ep_{\alpha}}{\partial x_3}=\frac{2}{\ep}e_{\alpha 3}(u^\ep)
\end{equation} 
is bounded in $L^2(\Omb)$, as a consequence of (\ref{97}), and, as 
$v^\ep_{\alpha}$ has meanvalue zero with respect to $x_3$, it results that it is bounded in 
$L^2(\omb;H^1(-1,0))$, so that (\ref{911}) holds true, i.e.
 \begin{equation}\label{936}
v^\ep_{\alpha} \rightharpoonup v^b_\alpha \mbox{ weakly in }L^2(\omb;H^1(-1,0)),
 \end{equation}
for some subsequence of $\ep$ and for some $v^b_\alpha$ in $L^2(\omb;H^1(-1,0))$. 
Setting $v^b=(v^b_1,v^b_2,0)$, we get 
\[
e_{\alpha 3}(v^b)=\frac{1}{2}\frac{\partial v^b_{\alpha}}{\partial x_3},
\]
so that we derive (\ref{98}) from (\ref{935}) and (\ref{936}). 

Finally we prove (\ref{99}) and (\ref{910}), by introducing the sequence of functions
\[
w^\ep=\frac{1}{\ep^2}u^\ep_3 -\int_{-1}^0\frac{1}{\ep^2}u^\ep_3 \, dx_3,
\]
which is bounded in $L^2(\omb;H^1(-1,0))$, since 
\[
\frac{\partial w^\ep}{\partial x_3} =\frac{1}{\ep^2} e_{33}(u^\ep)
\]
is  bounded in $L^2(\Omb)$, due to (\ref{97}). So, extracting a subsequence, 
we can find $w^b_3$ in $L^2(\omb;H^1(-1,0))$, having meanvalue zero in $x_3$, 
such that (\ref{99}) and (\ref{910}) hold true, which ends the proof of Lemma 8.

\subsection{The density arguments}

In Section 6, we have mentionned four density arguments. 
These are stated in the following lemmata and proved hereafter. This is done 
for sake of completeness, since Lemmata 7 and 8 are very classical, 
Lemma 5 is classical and very similar to the density result proved in 
\cite{GaGuLeMo2}, while Lemma 6, though less classical, results from 
Theorem 9.1.3 of \cite{AdHe}.

\begin{lemma}
Let $v \in H^1_0(\omb)$, $0 \in \omb \subset \RR^2$.  There exist a sequence 
of positive numbers $r^n$, tending to zero, and a sequence 
of functions $v^n \in H^1_0(\omb)$, such that
$$
v^n \equiv 0 \mbox{ in the ball } B^n \mbox{ of center $0$ and radius $r^n$ },
$$
$$
v^n \rightarrow v \mbox{ in } H^1_0(\omb).
$$
\end{lemma}

Proof: Let $\tilde{V}=\{v \in {\cal C}^1(\overline{\omb}), \, v=0 
\mbox{ on }\partial \omb \}$. As $\tilde{V}$ is dense in $H^1_0(\omb)$, 
we may restrict to $v$ in $\tilde{V}$. Then the proof goes as follows. For any 
integer $n$, we consider two balls $B^n$ and ${B'}^n$ in $\omb \subset \R^2$, 
with center $0$ and respective radii $r^n$ and $R^n$, to be determined 
later on, and such that $0<r^n<R^n$, $R^n$ tends to zero as $n$ tends to 
infinity. We define $v^n \in H^1_0(\omb)$ by:
\[
v^n =0 \mbox{ in }B^n, \; v^n =v\mbox{ in }\omb \setminus {B'}^n,\; 
v^n = (1-\phi^n)v\mbox{ in } {B'}^n \setminus B^n,
\]
where $\phi^n$ is the solution of the capacity problem in 
${B'}^n \setminus B^n$:
\[ 
\Delta \phi^n =0 \mbox{ in }{B'}^n \setminus B^n, 
\; \phi^n =1 \mbox{ on } \partial B^n,
\; \phi^n =0 \mbox{ on } \partial B'^n.
\]
It is clear that 
$v^n \in W^{1,\infty}(\omb) \cap H^1_0(\omb)$ and $v^n \equiv 0$ in $B^n$. We 
are going to prove that, for convenient $r^n$ and $R^n$, $v^n \rightarrow v$ in 
$H^1_0(\omb)$. Actually, as $0 \leq \phi^n \leq 1$, 
\[ \begin{array}{ll}
\|v^n-v\|^2_{H^1_0(\omb)}
&=\int_{{B'}^n}|\nabla (v^n-v)|^2\, dx'\\ \\
&=\int_{B^n}|\nabla v|^2\, dx'+\int_{{B'}^n \setminus B^n}
| \phi^n\nabla v +v\nabla \phi^n|^2\, dx'\\ \\
&\leq \int_{B^n}|\nabla v|^2\, dx'+ 2 \int_{{B'}^n \setminus B^n}
| \phi^n\nabla v|^2\, dx' +2 \int_{{B'}^n \setminus B^n}
|v\nabla \phi^n|^2\, dx'\\ \\
&\leq 3\int_{{B'}^n}|\nabla v|^2\, dx'+2 \int_{{B'}^n \setminus B^n}
|v\nabla \phi^n|^2\, dx'\\ \\
&\leq 3\pi R^n \|\nabla v\|^2_{L^\infty(\omb)} +2 \| v\|^2_{L^\infty(\omb)}
\int_{{B'}^n \setminus B^n}
|\nabla \phi^n|^2\, dx'\\ \\
&=3\pi R^n \|\nabla v\|^2_{L^\infty(\omb)} +4 \pi \| v\|^2_{L^\infty(\omb)} 
\left(\log \frac{R^n}{r^n} \right)^{-1}.
\end{array}
\]
It is enough to take (e.g.) $r^n=1/n^2$ and $R^n=1/n$, in order to get 
$v^n \rightarrow v$ in $H^1_0(\omb)$.

\begin{lemma}
Let $v \in H^2_0(\omb)$, $0 \in \omb \subset \RR^2$, $v(0)=0$.  There exist 
a sequence of positive numbers $r^n$, tending to zero, and a sequence 
of functions $v^n \in H^2_0(\omb)$, such that
$$
v^n \equiv 0 \mbox{ in the ball } B^n \mbox{ of center $0$ and radius $r^n$ },
$$
$$
v^n \rightharpoonup v \mbox{ weakly in } H^2_0(\omb).
$$
\end{lemma}

Proof: 1) For any $v\in H^2_0(\omb)$, with $v(0)=0$, there exists 
$\overline{v}^n\in{\cal C}^2(\overline{\omb})\cap H^2_0(\omb)$, such that 
$\overline{v}^n \rightarrow v$ in $H^2(\omb)$ and hence in 
${\cal C}^0(\overline{\omb})$. In particular, 
$\overline{v}^n(0)\rightarrow v(0)=0$. Setting 
$v^n=\overline{v}^n-\overline{v}^n(0)\phi$, with $\phi\in {\cal D}(\omb)$ and 
$\phi(0)=1$, it is clear that 
$v^n \in {\cal C}^2(\overline{\omb})\cap H^2_0(\omb)$, $v^n(0)=0$ and 
$v^n\rightarrow v$ in $H^2(\omb)$.

2) From step 1), we may restrict to $v$ in 
${\cal C}^2(\overline{\omb})\cap H^2_0(\omb)$, $v(0)=0$. Let $v^n=v\phi^n$, 
with $\phi^n(x')=\phi(n|x'|)$ and $\phi \in {\cal C}^\infty(\R)$, 
$0 \leq \phi \leq 1$, $\phi \equiv 0$ on $(-\infty,1]$, $\phi \equiv 1$ on 
$[2,+\infty)$. Clearly $v^n \in H^2_0(\omb)$, $v^n \equiv 0$ in the ball 
of center $0$ and radius $1/n$ and
\[
\int_{\omb}|v^n-v|^2\,dx'\leq \int_{|x'|<\frac{2}{n}}|v|^2\,dx'\rightarrow 0,
\]
that is $v^n\rightarrow v$ in $L^2(\omb)$. Hence the Lemma is proved, 
as soon as we have proved that $v^n$ is bounded uniformly in $H^2_0(\omb)$, i.e. 
\begin{equation}\label{d2}
\frac{\partial ^2 v^n}{\partial x_{\alpha}\partial x_{\beta}} 
\mbox{ is bounded in }L^2(\omb).
\end{equation}
But 
\[  
\frac{\partial ^2 v^n}{\partial x_{\alpha}\partial x_{\beta}} 
=v\frac{\partial ^2 \phi^n}{\partial x_{\alpha}\partial x_{\beta}}+
\phi^n\frac{\partial ^2 v}{\partial x_{\alpha}\partial x_{\beta}}+
\frac{\partial v}{\partial x_{\alpha}}
\frac{\partial\phi^n }{\partial x_{\beta}}+
\frac{\partial v}{\partial x_{\beta }}
\frac{\partial\phi^n }{\partial x_{\alpha}}.
\]
The second term is obviously bounded in $L^\infty(\omb)$. Moreover, since
\[ 
\frac{\partial\phi^n }{\partial x_{\alpha}}
=n\phi '(n|x'|)\frac{x_{\alpha}}{|x'|} \; \mbox{ and }\;
\frac{\partial ^2 \phi^n}{\partial x_{\alpha}\partial x_{\beta}}=
n^2\phi "(n|x'|)\frac{x_{\alpha}x_{\beta}}{|x'|^2}
+n\phi '(n|x'|)\left(\frac{\delta_{\alpha \beta}}{|x'|} - 
 \frac{x_{\alpha}x_{\beta}}{|x'|^3} \right),
\]
it follows that 
\[
\left|\frac{\partial\phi^n }{\partial x_{\alpha}}\right| \leq Cn \mbox{ and }
\left|\frac{\partial ^2 \phi^n}{\partial x_{\alpha}\partial x_{\beta}}\right|\leq Cn^2,
\]
\[
\int_{\omb} \left|\frac{\partial v}{\partial x_{\beta }}
\frac{\partial\phi^n }{\partial x_{\alpha}}\right|^2 \, dx'\leq 
C  \left\|\frac{\partial v}{\partial x_{\beta }}\right\|^2_\infty
\int_{\frac{1}{n}<|x'|<\frac{2}{n}}n^2\, dx'=
C  \left\|\frac{\partial v}{\partial x_{\beta }}\right\|^2_\infty
\int_{ 1<|x'|< 2} \, dx'=C,
\]
\[
\int_{\omb}\displaystyle 
|v\frac{\partial ^2 \phi^n}{\partial x_{\alpha}\partial x_{\beta}}|^2 \, dx'
\leq \|v \|^2_{L^\infty(\frac{1}{n}<|x'|<\frac{2}{n})}
\int_{\frac{1}{n}<|x'|<\frac{2}{n}} Cn^4\, dx'=Cn^2
\|v \|^2_{L^\infty(\frac{1}{n}<|x'|<\frac{2}{n})}.
\]
But, for $1/n<|x'|<2/n$, $|v(x')| \leq C |x'|\leq C/n$, since 
$v$ is regular and $v(0)=0$. It follows that 
\[
\int_{\omb}\displaystyle 
|v\frac{\partial ^2 \phi^n}{\partial x_{\alpha}\partial x_{\beta}}|^2 \, dx'
 \leq C
\]
and finally, (\ref{d2}) holds true, ending the proof of Lemma 5.

\begin{lemma}
Let $v \in L^2(\omb;H^1(-1,0))$, $0 \in \omb \subset \RR^2$.  There exist 
a sequence of positive numbers $r^n$, tending to zero, and a sequence 
of functions $v^n$, such that
$$
v^n \in {\cal C}^1(\overline{\Omb}),
$$
$$
v^n \equiv 0 \mbox{ in } B^n \times \{0\}, B^n \mbox{ denoting the ball of 
center $0$ and radius $r^n$ },
$$
$$
v^n \rightarrow v \mbox{ in } L^2(\omb;H^1(-1,0)).
$$
\end{lemma}

Proof: By density of ${\cal C}^1(\overline{\Omb})$ in $L^2(\omb;H^1(-1,0))$, 
we may restrict to $v \in {\cal C}^1(\overline{\Omb})$. We consider a 
sequence $r^n$ of positive numbers, converging to zero, 
and a sequence of 
functions $\phi^n:\omb \rightarrow \R$, of class ${\cal C}^\infty$, with 
$\phi^n \equiv 0$ in the ball $B^n$ of center $0$ and radius $r^n$, 
$\phi^n \equiv 1$ outside the ball $B'^n$ of center $0$ and radius $2 r^n$, 
$0\leq \phi^n \leq1$ in $B'^n \setminus B^n$. We set $v^n=\phi^n v$. Then clearly 
$v^n \in {\cal C}^1(\overline{\Omb})$ and
\[\begin{array}{ll}  
\|v^n -v\|^2_{L^2(\omb;H^1(-1,0))}&=\int_{\Omb}|v^n -v|^2 \, dx +
\int_{\Omb}|\frac{\partial}{dx_3}(v^n -v)|^2 \, dx \\ \\
&= \int_{B^n \times (-1,0)}|v|^2 \, dx +  
\int_{(B'^n \setminus B^n) \times (-1,0)}|(1-\phi^n) v|^2 \, dx +\\ \\
&+ \int_{B^n \times (-1,0)}|\frac{\partial v}{\partial x_3}|^2 \, dx +   
\int_{(B'^n \setminus B^n) \times (-1,0)}|(1-\phi^n) 
\frac{\partial v}{\partial x_3}|^2 \, dx \\ \\
&\leq \int_{B'^n \times (-1,0)} \left( |v|^2 + 
|\frac{\partial v}{\partial x_3}|^2\right)\, dx,
\end{array}
\]
which tends to zero, as soon as $r^n$ tends to zero.
 
\begin{lemma}
Let $U=H^1_T(0,1)=\{u \in H^1(0,1), \;u(1)=0\}$, 
$\tilde{U}=\{u \in {\cal C}^1[0,1], \;u(1)=0\}$, 
$V=H^2_0(\omb)$, $0 \in \omb \subset \RR^2$, $\tilde{V}={\cal C}^1(\overline{\omb}) \cap H^2_0(\omb)$. Now 
let $W=\{(u,v) \in U \times V, \; u(0)=v(0)\}$, 
$\tilde{W}=\{(u,v) \in \tilde{U} \times \tilde{V}, \; u(0)=v(0)\}$. 
Then $\tilde{W}$ is dense in $W$.
\end{lemma}

Proof: It is clear that $\tilde{U}$ is dense in $U$ and that $\tilde{V}$ 
is dense in $V$. Therefore, for any $(u,v) \in W$, there exists 
$(\overline{u}^n,\overline{v}^n) \in \tilde{U} \times \tilde{V}$ such that 
\[\begin{array}{l}
\overline{u}^n \rightarrow u \mbox{ in }H^1(0,1)
\mbox{ and hence in } {\cal C}^0[0,1],\\
\overline{v}^n \rightarrow v \mbox{ in }H^2(\omb)
\mbox{ and hence in } {\cal C}^0(\overline{\omb}).
\end{array}
\]
Let $\phi^1\in {\cal C}^\infty[0,1]$ with $\phi^1(0)=1$, $\phi^1(1)=0$, 
$\phi^2\in {\cal D}(\omb)$ with $\phi^2(0)=1$ and let 
\[
u^n=\overline{u}^n -(\overline{u}^n(0)-u(0))\phi^1,
\]
\[
v^n=\overline{v}^n -(\overline{v}^n(0)-v(0))\phi^2.
\]
It is clear that $u^n\in \tilde{U}$, $v^n\in \tilde{V}$ and 
$u^n(0)=u(0)=v(0)=v^n(0)$, so that $(u^n,v^n)\in \tilde{W}$. Moreover
\[
\|u^n-\overline{u}^n\|_{H^1(0,1)}=|\overline{u}^n(0)-u(0)|\|\phi^1 \|_{H^1(0,1)} 
\rightarrow 0
\]
and hence $u^n\rightarrow u$ in $H^1(0,1)$. Similarly $v^n\rightarrow v$ 
in $H^2(\omb)$.

{\bf Acknowledgment:} This work was partially supported by the 2001 
and 2002 Vinci programs.

\end{document}